\def\DateTime{22/December/1997, 7:45PM}
\def\Version{Version 4.5}
\def\yes{\if00}

\def\ifpicture{\yes}
\def\iftwelvept{\yes}

\iftwelvept
\documentclass[leqno,12pt]{amsart}
\else
\documentclass[leqno]{amsart}
\fi
\usepackage{amssymb}
\usepackage{amscd}
\usepackage{latexsym}

\iftwelvept
\else\fi

\theoremstyle{plain}
\newtheorem{Theorem}{Theorem}[section]
\newtheorem{Proposition}[Theorem]{Proposition}
\newtheorem{Lemma}[Theorem]{Lemma}
\newtheorem{Corollary}[Theorem]{Corollary}
\newtheorem{Claim}{Claim}[Theorem]

\theoremstyle{definition}

\newtheorem{Remark}[Theorem]{Remark}
\newtheorem{Example}[Theorem]{Example}
\newtheorem{Conjecture}[Theorem]{Conjecture}

\renewcommand{\theTheorem}{\arabic{section}.\arabic{Theorem}}
\renewcommand{\theClaim}{\arabic{section}.\arabic{Theorem}.\arabic{Claim}}
\renewcommand{\theequation}{\arabic{section}.\arabic{Theorem}.\arabic{Claim}}

\def\rom{\textup}

\newcommand{\QQ}{{\mathbb{Q}}}
\newcommand{\RR}{{\mathbb{R}}}

\newcommand{\PP}{{\mathbb{P}}}
\newcommand{\OO}{{\mathcal{O}}}
\newcommand{\rank}{\operatorname{rk}}
\newcommand{\codim}{\operatorname{codim}}
\newcommand{\ch}{\operatorname{char}}
\newcommand{\dis}{\operatorname{dis}}
\newcommand{\Ker}{\operatorname{Ker}}

\newcommand{\Jac}{\operatorname{Jac}}
\newcommand{\Sym}{\operatorname{Sym}}
\newcommand{\Hilb}{\operatorname{Hilb}}
\newcommand{\Supp}{\operatorname{Supp}}

\newcommand{\End}{\operatorname{\mathcal{E}\textsl{nd}}}
\newcommand{\Pic}{\operatorname{Pic}}
\newcommand{\Div}{\operatorname{Div}}
\newcommand{\Vt}{\operatorname{Vert}}
\newcommand{\Ed}{\operatorname{Ed}}
\newcommand{\Sing}{\operatorname{Sing}}
\newcommand{\Spec}{\operatorname{Spec}}
\newcommand{\Chow}{A}
\newcommand{\chern}{\operatorname{ch}}
\newcommand{\NS}{NS}
\newcommand{\Nef}{\operatorname{Nef}}
\newcommand{\WP}{\operatorname{WP}}
\newcommand{\Curve}{\operatorname{Cur}}
\newcommand{\Proj}{\operatorname{Proj}}
\newcommand{\zero}{\operatorname{div}}
\newcommand{\Proof}{{\sl Proof.}\quad}
\newcommand{\QED}{{\unskip\nobreak\hfil\penalty50\quad\null\nobreak\hfil
{$\Box$}\parfillskip0pt\finalhyphendemerits0\par\medskip}}
\newcommand{\rest}[2]{\left.{#1}\right\vert_{{#2}}}

\begin{document}

\title[Relative Bogomolov's inequality]%
{Relative Bogomolov's inequality \\
and the cone of positive divisors \\
on the moduli space of stable curves}
\author{Atsushi Moriwaki}
\address{Department of Mathematics, Faculty of Science,
Kyoto University, Kyoto, 606-01, Japan}
\email{moriwaki@kusm.kyoto-u.ac.jp}
\keywords{Bogomolov's inequality, moduli space, stable curve}
\subjclass{Primary 14H10, 14C20; Secondary 14G40}
\date{\DateTime, (\Version)}
\thanks{This paper will appear in Journal of AMS}
\begin{abstract}
Let $f : X \to Y$ be a surjective and
projective morphism of smooth quasi-projective varieties
over an algebraically closed field of characteristic zero
with $\dim f = 1$. 
Let $E$ be a vector bundle of rank $r$ on $X$. 
In this paper, we would like to show that
if $X_y$ is smooth and $E_y$ is semistable for some $y \in Y$, then
$f_*\left( 2rc_2(E) - (r-1)c_1(E)^2 \right)$
is weakly positive at $y$.
We apply this result to obtain the following description of the
cone of weakly positive $\QQ$-Cartier divisors
on the moduli space of stable curves.
Let $\overline{\mathcal{M}}_g$ (resp. $\mathcal{M}_g$) be the moduli space of 
stable (resp. smooth) curves of genus $g \geq 2$.
Let $\lambda$ be the Hodge class and
$\delta_i$'s ($i = 0, \ldots, [g/2]$) the boundary classes.
Then, a $\QQ$-Cartier divisor $x \lambda + \sum_{i=0}^{[g/2]} y_i \delta_i$
on $\overline{\mathcal{M}}_g$ is weakly positive over $\mathcal{M}_g$ if and only if
$x \geq 0$, $g x + (8g + 4) y_0 \geq 0$, and
$i(g-i) x + (2g+1) y_i \geq 0$ for all $1 \leq i \leq [g/2]$.
\end{abstract} 

\maketitle


\section*{Introduction}
\renewcommand{\theTheorem}{\Alph{Theorem}}
Throughout this paper, we fix an algebraically closed field $k$.

Let $f : X \to Y$ be a surjective and
projective morphism of quasi-projective varieties over $k$ with $\dim f = 1$.
Let $E$ be a vector bundle of rank $r$ on $X$. Then, we define
the {\em discriminant divisor} of $E$
with respect to $f : X \to Y$ to be
\[
\dis_{X/Y}(E) = f_*\left( \left(2rc_2(E) - (r-1)c_1(E)^2 \right) \cap [X] \right).
\]
Here $f_*$ is the push-forward of cycles, so that
$\dis_{X/Y}(E)$ is a divisor modulo linear equivalence on $Y$.
In this paper, we would like to show the following theorem
(cf. Corollary~\ref{cor:nef:psudo:dis}) and give its applications.

\begin{Theorem}[$\ch(k) = 0$]
\label{thm:intro:A}
We assume that $Y$ is smooth over $k$.
Let $y$ be a point of $Y$ and $\overline{\kappa(y)}$ the
algebraic closure of the residue field $\kappa(y)$ at $y$.
If $f$ is flat over $y$, 
the geometric fiber $X_{\bar{y}} = X \times_Y \Spec\left(\overline{\kappa(y)}\right)$
over $y$ is reduced and Gorenstein,
and $E$ is semistable on each
connected component of the normalization of $X_{\bar{y}}$, 
then $\dis_{X/Y}(E)$ is weakly positive at $y$, namely,
for any ample divisors $A$ on $Y$ and any positive integers $n$,
there is a positive integer $m$ such that
\[
H^0(Y, \OO_Y(m(n \dis_{X/Y}(E) + A))) \otimes \OO_Y \to
\OO_Y(m(n \dis_{X/Y}(E) + A))
\]
is surjective at $y$.
Note that this theorem still holds in positive characteristic under
the strong semistability of $E_{\bar{y}}$ 
\rom{(}cf. Corollary~\rom{\ref{cor:nef:psudo:dis:in:p}}\rom{)}.
\end{Theorem}

An interesting point of the above theorem is that
even if the weak positivity of $\dis_{X/Y}(E)$ at $y$ 
is a global property on $Y$,
it can be derived from the local assumption
``the goodness of $X_{\bar{y}}$ and the semistability of $E_{\bar{y}}$''.
This gives a great advantage to our applications.

In order to understand the intuition underlying the theorem, 
let us consider a toy case.
Namely, we suppose that $f : X \to Y$ is a smooth surface fibred  over a curve
and the fiber is general.
Bogomolov's instability theorem \cite{Bogo} says that
if $E_{\bar{y}}$ is semistable, then the codimension two cycle
$2r c_2(E) - (r-1)c_1(E)^2$ has non-negative degree.
So if we push it down to a codimension one cycle on $Y$, then
one can rephrase Bogomolov's theorem as saying that
the semistability of $E_y$ implies the non-negativity of $\dis_{X/Y}(E)$.

\medskip
An immediate application of our inequality
is a solution concerning the positivity of
divisors on the moduli space of stable curves.
Let $g \geq 2$ be an integer, and
$\overline{\mathcal{M}}_g$ (resp. $\mathcal{M}_g$) the moduli space of 
stable (resp. smooth) curves of genus $g$ over $k$.
The boundary $\overline{\mathcal{M}}_g \setminus \mathcal{M}_g$
is of codimension one and has $[g/2]+1$ irreducible components, say,
$\Delta_0, \Delta_1, \ldots, \Delta_{[g/2]}$.
The geometrical meaning of indexes is as follows.
A general point of $\Delta_0$ represents
an irreducible stable curve with one node, and
a general point of $\Delta_i$ ($i > 0$) represents
a stable curve consisting of a curve of genus $i$ and
a curve of genus $g-i$ joined at one point.
Let $\delta_i$ be the class of $\Delta_i$ in
$\Pic(\overline{\mathcal{M}}_g) \otimes \QQ$
(strictly speaking, $\delta_i = c_1(\OO(\Delta_i))$ for $i \not= 1$, and
$\delta_1 = \frac{1}{2} c_1(\OO(\Delta_1))$),
and $\lambda$ the Hodge class on 
$\overline{\mathcal{M}}_g$.
A fundamental problem due to Mumford \cite{Mum} is to decide which $\QQ$-divisor
\[
a \lambda - b_0 \delta_0 - b_1 \delta_1 - \cdots - b_{[g/2]}\delta_{[g/2]}
\]
is positive, where $a, b_0, \ldots, b_{[g/2]}$ are rational numbers. 
Here, we can use a lot of types of positivity,
namely, ampleness, numerical effectivity, effectivity, pseudo-effectivity, and so on.
Besides them, we would like to introduce a new sort of positivity for our purposes. 
Let $V$ be a projective variety over $k$ and $U$ a
non-empty Zariski open set of $V$.
A $\QQ$-Cartier divisor $D$ on $V$ is said to be {\em numerically effective over $U$}
if $(D \cdot C) \geq 0$ for all irreducible curves $C$ on $V$ with 
$C \cap U \not= \emptyset$.
A first general result in this direction was found by
Cornalba-Harris \cite{CH}, Xiao \cite{Xi} and Bost \cite{Bo}.
They proved that the $\QQ$-divisor
\[
(8g+4) \lambda - g (\delta_0 + \delta_1 + \cdots + \delta_{[g/2]})
\]
is numerically effective over $\mathcal{M}_g$.
As we observed in \cite{Mo4} and \cite{Mo5}, 
it is not sharp in coefficients of $\delta_i$ ($i > 0$).
Actually, the existence of a certain refinement of the above result
was predicted at the end of the paper \cite{CH}.
Our solution for this problem is the following
(cf. Theorem~\ref{thm:wpos:at:M:g} and Proposition~\ref{prop:samp:wp:U}).

\begin{Theorem}[$\ch(k) = 0$]
\label{thm:intro:wp}
The divisor
\[
(8g+4) \lambda - g \delta_0 - \sum_{i=1}^{[g/2]} 4 i (g-i) \delta_i
\]
is weakly positive over $\mathcal{M}_g$, i.e.,
if we denote the above divisor by $D$, 
then for any ample $\QQ$-Cartier divisors $A$ on $\overline{\mathcal{M}}_g$,
there is a positive integer $n$ such that
$n(D + A)$ is a Cartier divisor and
\[
H^0(\overline{\mathcal{M}}_g, \OO_{\overline{\mathcal{M}}_g}(n(D+A))) \otimes
\OO_{\overline{\mathcal{M}}_g} \to \OO_{\overline{\mathcal{M}}_g}(n(D+A))
\]
is surjective on $\mathcal{M}_g$.
In particular, it
is pseudo-effective, and numerically effective over $\mathcal{M}_g$.
\end{Theorem}

As an application of this theorem,
we can decide the cone of weakly positive divisors over $\mathcal{M}_g$
(cf. Corollary~\ref{cor:wp:cone:equal}).

\begin{Theorem}[$\ch(k) = 0$]
If we denote by
$\WP(\overline{\mathcal{M}}_g; \mathcal{M}_g)$
the cone in $\Pic(\overline{\mathcal{M}}_g) \otimes \QQ$
consisting of weakly positive $\QQ$-Cartier divisors over $\mathcal{M}_g$, then
\[
\WP(\overline{\mathcal{M}}_g; \mathcal{M}_g)
= \left\{ x \lambda + \sum_{i=0}^{[g/2]} y_i \delta_i \ \left| \ 
\begin{array}{l}
x \geq 0, \\
g x + (8g + 4) y_0 \geq 0, \\
i(g-i) x + (2g+1) y_i \geq 0 \quad (1 \leq i \leq [g/2]).
\end{array}
\right. \right\}.
\]
\end{Theorem}

Moreover, using Theorem~\ref{thm:intro:wp}, 
we can deduce a certain kind of inequality on an algebraic surface.
In order to give an exact statement, we will introduce types of nodes
of semistable curves.
Let $Z$ be a semistable curve over $k$, and
$P$ a node of $Z$. 
We can assign a number $i$ to the node $P$ in the following way.
Let $\iota_P : Z_P \to Z$ be the partial normalization of $Z$
at $P$. If $Z_P$ is connected, then $i=0$. 
Otherwise, $i$ is the minimum of arithmetic genera of 
two connected components of $Z_P$. We say the node $P$ of
$Z$ is {\em of type $i$}.

Let $X$ be a smooth projective surface over $k$, 
$Y$ a smooth projective curve over $k$, and
$f : X \to Y$ a semistable curve of genus $g \geq 2$ over $Y$.
By abuse of notation, we denote by $\delta_i(X/Y)$ the number of
nodes of type $i$ in all singular fibers of $f$.
Actually, $\delta_i(X/Y) = \deg(\pi^*(\delta_i))$, where
$\pi : Y \to \overline{\mathcal{M}}_g$ is the morphism induced by $f : X \to Y$.
Then, we have the following
(cf. Corollary~\ref{cor:sharp:slope:inq}).

\begin{Theorem}[$\ch(k) = 0$]
\label{thm:sharp:slope:inq:in:intro}
With notation being as above, we have the inequality
\[
(8g+4) \deg(f_*(\omega_{X/Y})) \geq
g \delta_0(X/Y) + \sum_{i=1}^{[g/2]} 4 i (g-i) \delta_i(X/Y).
\]
\end{Theorem}

As an arithmetic application of
Theorem~\ref{thm:sharp:slope:inq:in:intro}, 
we can show the following answer
for effective Bogomolov's conjecture over function fields
(cf. Theorem~\ref{thm:bogomolov:function:field}).
(Recently, Bogomolov's conjecture over number fields was
solved by Ullmo \cite{Ul}, but effective Bogomolov's conjecture
is still open.)

\begin{Theorem}[$\ch(k) = 0$]
We assume that $f$ is not smooth and
every singular fiber of $f$
is a tree of stable components, i.e.,
every node of type $0$ on the stable model of each singular fiber is
a singularity of an irreducible component, then 
effective Bogomolov's conjecture holds for the generic fiber of $f$. 
Namely, 
let $K$ be the function field of $Y$, $C$ the generic fiber of $f$,
$\Jac(C)$ the Jacobian of $C$, and let $j : C(\overline{K})
\to \Jac(C)(\overline{K})$ be the morphism given by
$j(x) = (2g-2)x - \omega_C$. 
Then, the set $\{ x \in C(\overline{K}) \mid
\Vert j(x) - P \Vert_{NT} \leq r \}$ is finite for
any $P \in \Jac(C)(\overline{K})$ and any non-negative real numbers
$r$ less than
\[
\sqrt{\frac{(g-1)^2}{g(2g+1)}\left(
\frac{g-1}{3}\delta_0(X/Y) + \sum_{i=1}^{\left[\frac{g}{2}\right]}
4i(g-i)\delta_i(X/Y) \right)},
\]
where $\Vert \ \Vert_{NT}$ is 
the semi-norm arising from the Neron-Tate height paring
on $\Jac(C)(\overline{K})$.
\end{Theorem}

Finally, we would like to express our hearty thanks to
Institut des Hautes \'{E}tudes Scientifiques where all works of this paper
had been done, and to Prof. Bost who pointed out a fatal error of the previous
version of it.
We are also grateful to referees for their wonderful suggestions.

\section{Elementary properties of semi-ampleness and weak positivity}
\label{sec:pef:pamp:div}
\renewcommand{\theTheorem}{\arabic{section}.\arabic{Theorem}}
In this section,
we will introduce two kinds of positivity of divisors, namely
semi-ampleness and weak positivity, and investigate their
elementary properties.

Let $X$ be a $d$-dimensional algebraic variety over $k$.
Let $Z_{d - 1}(X)$ be a free abelian group generated by
integral subvarieties of dimension $d - 1$, and
$\Div(X)$ a group consisting of
Cartier divisors on $X$.
We denote $Z_{d - 1}(X)$ (resp. $\Div(X)$) modulo linear equivalence by
$\Chow_{d - 1}(X)$ (resp. $\Pic(X)$).
An element of $Z_{d - 1}(X) \otimes \QQ$ (resp. $\Div(X) \otimes \QQ$)
is called a {\em $\QQ$-divisor} (resp. {\em $\QQ$-Cartier divisor}) on $X$.

We say a $\QQ$-Cartier divisor $D$ is {\em the limit of 
a sequence $\{ D_m \}_{m=1}^{\infty}$ of $\QQ$-Cartier divisors in $\Pic(X) \otimes \QQ$},
denoted by ${\displaystyle D = \lim_{m \to \infty} D_m}$ in $\Pic(X) \otimes \QQ$,
if there are $\QQ$-Cartier divisors $Z_1, \ldots, Z_{l}$ and
infinite sequences $\{ a_{1, m} \}_{m=1}^{\infty}, \ldots, \{a_{l, m} \}_{m=1}^{\infty}$
of rational numbers such that (1) $l$ does not depend on $m$,
(2) $D = D_m + \sum_{i=1}^{l} a_{i, m} Z_i$ in $\Pic(X) \otimes \QQ$
for all $m \geq 1$, and
(3) ${\displaystyle \lim_{m\to\infty} a_{i, m} = 0}$ for all $i=1, \ldots, l$.
For example, a pseudo-effective $\QQ$-Cartier divisor is the limit of
effective $\QQ$-Cartier divisors in $\Pic(X) \otimes \QQ$.

Let $x$ be a point of $X$.
A $\QQ$-Cartier divisor $D$ on $X$ is said to be {\em semi-ample at $x$}
if there is a positive integer $n$ such that $nD \in \Div(X)$ and
$H^0(X, \OO_X(nD)) \otimes \OO_X \to \OO_X(nD)$ is surjective at $x$.
Further, according to Viehweg, $D$ is said to be {\em weakly positive at $x$}
if there is an infinite sequence $\{ D_m \}_{m=1}^{\infty}$ of
$\QQ$-Cartier divisors on $X$ such that
$D_m$ is semi-ample at $x$ for all $m \geq 1$ 
and ${\displaystyle D = \lim_{m \to \infty} D_m}$ in $\Pic(X) \otimes \QQ$. 
It is easy to see that
if $D$ is weakly positive at $x$, then $(D \cdot C) \geq 0$
for any complete irreducible curves $C$ passing through $x$.
As compared with the last property,
weak positivity has an advantage that
we can avoid bad subvarieties of codimension two
(cf. Proposition~\ref{prop:wp:codim:2}).

In order to consider properties of semi-ample or weakly positive divisors,
let us begin with the following two lemma.

\begin{Lemma}[$\ch(k) \geq 0$]
\label{lem:gen:by:global:sec:finite}
Let $\pi : X \to Y$ be a proper morphism of quasi-projective varieties over $k$ and
$y$ a point of $Y$ such that $\pi$ is finite over $y$.
Let $F$ be a coherent $\OO_X$-module and $H$ an ample line bundle on $Y$.
Then there is a positive integer $n_0$ such that,
for all $n \geq n_0$,
\[
 H^0(X, F \otimes \pi^{*}(H^{\otimes n})) \otimes \OO_X \to F \otimes \pi^{*}(H^{\otimes n})
\]
is surjective at each point of $\pi^{-1}(y)$.
\end{Lemma}

\Proof
Let $n_0$ be a positive integer such that,
for all $n \geq n_0$,
$\pi_*(F) \otimes H^{\otimes n}$ is generated by global sections, i.e.,
\[
 H^0(Y, \pi_*(F) \otimes H^{\otimes n}) \otimes \OO_Y \to \pi_*(F) \otimes H^{\otimes n}
\]
is surjective.
Thus,
\[
 H^0(X, F \otimes \pi^{*}(H^{\otimes n})) \otimes \OO_X \to \pi^* \pi_*(F \otimes \pi^*(H^{\otimes n}))
\]
is surjective because $\pi_*(F \otimes \pi^{*}(H^{\otimes n})) = 
\pi_*(F) \otimes H^{\otimes n}$. On the other hand, since $\pi$ is finite over $y$,
\[
\pi^* \pi_*(F \otimes \pi^*(H^{\otimes n}))
\to
F \otimes \pi^*(H^{\otimes n})
\]
is surjective at each point of $\pi^{-1}(y)$. 
Thus, we get our assertion.
\QED

\begin{Lemma}[$\ch(k) \geq 0$]
\label{lem:pamp:plus:good:samp}
Let $\pi : X \to Y$ be a proper morphism of quasi-projective varieties over $k$ and
$x$ a point of $X$ such that $\pi$ is finite over $\pi(x)$.
Let $D$ be a $\QQ$-Cartier divisor on $X$ and
$A$ a $\QQ$-Cartier divisor on $Y$.
If $D$ is weakly positive at $x$ and $A$ is ample, then
$D + \pi^*(A)$ is semi-ample at $x$.
\end{Lemma}

\Proof
By our assumption, there are $\QQ$-Cartier divisors $Z_1, \ldots, Z_{l}$,
an infinite sequence $\{ D_m \}_{m=1}^{\infty}$ of $\QQ$-Cartier divisors, and
infinite sequences $\{ a_{1, m} \}_{m=1}^{\infty}, \ldots, \{ a_{l, m} \}_{m=1}^{\infty}$
of rational numbers with the following properties.
\begin{enumerate}
\renewcommand{\labelenumi}{(\roman{enumi})}
\item
$D = D_m + \sum_{i=1}^{l} a_{i, m} Z_i$ in $\Pic(X) \otimes \QQ$
for sufficiently large $m$.

\item
$\lim_{m\to\infty} a_{i, m} = 0$ for all $i=1, \ldots, l$.

\item
$D_m$ is semi-ample at $x$ for $m \gg 0$.
\end{enumerate}
By virtue of Lemma~\ref{lem:gen:by:global:sec:finite},
we can find a positive integer $n$ such that
$n \pi^*(A) + Z_i$ and $n \pi^*(A) - Z_i$ are semi-ample at $x$.
We choose a sufficiently large $m$ with $|nla_{i, m}| < 1$ ($i=1, \ldots, l$).
Then,
\[
\frac{1 + nla_{i, m}}{2ln} > 0
\quad\text{and}\quad
\frac{1 - nla_{i, m}}{2ln} > 0.
\]
On the other hand,
\[
D + \pi^*(A) \sim
D_m + \sum_{i=1}^l \left(
\frac{1 + nla_{i, m}}{2ln}(n\pi^*(A) + Z_i) +
\frac{1 - nla_{i, m}}{2ln}(n\pi^*(A) - Z_i)
\right).
\]
Thus, $D + \pi^*(A)$ is semi-ample at $x$.
\QED

As immediate consequences of Lemma~\ref{lem:pamp:plus:good:samp},
we have the following propositions.
The first one is a characterization of weak positivity
in terms of ample divisors.

\begin{Proposition}[$\ch(k) \geq 0$]
\label{prop:criterion:pamp}
Let $X$ be a quasi-projective variety over $k$, $x$ a
point of $X$, and
$D$ a $\QQ$-Cartier divisor on $X$.
Then, the following are equivalent.
\begin{enumerate}
\renewcommand{\labelenumi}{(\arabic{enumi})}
\item
$D$ is weakly positive at $x$.

\item
For any ample $\QQ$-Cartier divisors $A$ on $X$,
$D + A$ is semi-ample at $x$.

\item
There is an ample $\QQ$-Cartier divisor $A$ on $X$
such that 
$D + \epsilon A$ is semi-ample at $x$ for any positive rational numbers $\epsilon$.
\end{enumerate}
\end{Proposition}

\begin{Proposition}[$\ch(k) \geq 0$]
\label{prop:wp:codim:2}
Let $X$ be a normal quasi-projective variety over $k$,
$X_0$ a Zariski open set of $X$, and $x$ a point of $X_0$.
Let $D$ be a $\QQ$-Cartier divisor on $X$ and $D_0 = \rest{D}{X_0}$.
If $\codim(X \setminus X_0) \geq 2$, then we have the following.
\begin{enumerate}
\renewcommand{\labelenumi}{(\arabic{enumi})}
\item
$D$ is semi-ample at $x$ if and only if
$D_0$ is semi-ample at $x$.

\item
$D$ is weakly positive at $x$ if and only if
$D_0$ is weakly positive at $x$.
\end{enumerate}
\end{Proposition}

Next, let us consider functorial properties of semi-ampleness and weak positivity
under pull-back and push-forward.

\begin{Proposition}[$\ch(k) \geq 0$]
\label{prop:samp:pamp:pullback}
Let $\pi : X \to Y$ be a morphism of quasi-projective varieties.
Let $D$ be a $\QQ$-Cartier divisor on $Y$ and $x$ a point of $X$.
If $\pi^*(D)$ is defined, then we have the following.
\rom{(}Note that even if $\pi^*(D)$ is not defined, there is a $\QQ$-Cartier
divisor $D'$ such that $D' \sim D$ and $\pi^*(D')$ is defined.\rom{)}
\begin{enumerate}
\renewcommand{\labelenumi}{(\arabic{enumi})}
\item
If $D$ is semi-ample at $\pi(x)$, then
$\pi^{*}(D)$ is semi-ample at $x$.

\item
If $D$ is weakly positive at $\pi(x)$, then
$\pi^*(D)$ is weakly positive at $x$.
\end{enumerate}
\end{Proposition}

\Proof
(1) By our assumption,
$H^0(Y, \OO_Y(nD)) \otimes \OO_Y \to \OO_Y(nD)$ is surjective at $\pi(x)$
for a sufficiently large $n$.
Thus,
$H^0(Y, \OO_Y(nD)) \otimes \OO_X \to \OO_X(n\pi^*(D))$ is
surjective at $x$.
Here let us consider the following commutative diagram:
\[
\begin{CD}
H^0(Y, \OO_{Y}(nD)) \otimes \OO_{X} @>{\alpha}>> \OO_X(n\pi^*(D)) \\
@VVV @| \\
H^0(X, \OO_{X}(n\pi^*(D))) \otimes \OO_{X} @>{\alpha'}>> \OO_X(n\pi^*(D)).
\end{CD}
\]
Since $\alpha$ is surjective at $x$, so is $\alpha'$.
Therefore, $\pi^*(D)$ is semi-ample at $x$.

\medskip
(2) Let $A$ be an ample divisor on $Y$ such that
$\pi^*(A)$ is defined. Then, by Lemma~\ref{lem:pamp:plus:good:samp},
$D + (1/n)A$ is semi-ample at $\pi(x)$ for all $n > 0$.
Thus, by (1), $\pi^*(D) + (1/n)\pi^*(A)$ is semi-ample at $x$
for all $n > 0$. 
Therefore, $\pi^*(D)$ is weakly positive at $x$.
\QED

\begin{Proposition}[$\ch(k) \geq 0$]
\label{prop:samp:pamp:push}
Let $\pi : X \to Y$ be a surjective, proper and generically finite morphism
of normal quasi-projective varieties over $k$.
Let $D$ be a $\QQ$-Cartier divisor on $X$ and $y$ a point of $Y$ such that
$\pi_*(D)$ is a $\QQ$-Cartier divisor on $Y$ and
$\pi$ is finite over $y$. We set $\pi^{-1}(y) = \{x_1, \ldots, x_n \}$.
Then, we have the following.
\begin{enumerate}
\renewcommand{\labelenumi}{(\arabic{enumi})}
\item
If $D$ is semi-ample at $x_1, \ldots, x_n$, then
$\pi_*(D)$ is semi-ample at $y$.

\item
If $D$ is weakly positive at $x_1, \ldots, x_n$, then
$\pi_*(D)$ is weakly positive at $y$.
\end{enumerate}
\end{Proposition}

\Proof
(1) Clearly, we may assume that $D$ is a Cartier divisor.
If we take a sufficiently large integer $m$, then
$H^0(X, \OO_X(mD)) \otimes \OO_X \to \OO_X(mD)$
is surjective at $x_1, \ldots, x_n$.
Thus, there are sections $s_1, \ldots, s_n$ of $H^0(X, \OO_X(m D))$
with $s_i(x_i) \not= 0$ for all $i=1, \ldots, n$.
For $\alpha = (\alpha_1, \ldots\, \alpha_n) \in k^n$,
we set $s_{\alpha} = \alpha_1 s_1 + \cdots + \alpha_n s_n$.
Further, we set $V_i = \{ \alpha \in k^n \mid s_{\alpha}(x_i) = 0 \}$.
Then, $\dim V_i= n-1$ for all $i$. Thus, since
$\#(k) = \infty$, there is $\alpha \in k^n$
with $\alpha \not\in V_1 \cup \cdots \cup V_r$, i.e.,
$s_{\alpha}(x_i) \not= 0$ for all $i$.
Let us consider a divisor $E = \operatorname{div}(s_{\alpha})$. Then, $E \sim m D$.
Thus, $\pi_*(E) \sim m \pi_*(D)$. Here, $x_i \not\in E$
for all $i$. Hence, $y \not\in \pi_*(E)$. Therefore, we get our assertion.

\medskip
(2) Let $A$ be an ample divisor on $Y$. We set $D_m = D + (1/m)\pi^*(A)$.
Then, by Lemma~\ref{lem:pamp:plus:good:samp},
$D_m$ is semi-ample at $x_1, \ldots, x_n$.
Thus, by (1), $\pi_*(D_m) = \pi_*(D) + (1/m) \deg(\pi) A$ is semi-ample at $y$.
Therefore, $\pi_*(D)$ is weakly positive at $y$. 
\QED

\bigskip
Finally, let us consider semi-ampleness and weak positivity over an open set.
Let $X$ be a quasi-projective variety over $k$,
$U$ a Zariski open set of $X$, and $D$ a $\QQ$-Cartier divisor on $X$.
We say $D$ is {\em semi-ample over $U$} (resp. {\em weakly positive over $U$})
if $D$ is semi-ample (resp. weakly positive) at all points of $U$.
Then, we can easily see the following.

\begin{Proposition}[$\ch(k) \geq 0$]
\label{prop:samp:wp:U}
\begin{enumerate}
\renewcommand{\labelenumi}{(\arabic{enumi})}
\item
If $D$ is semi-ample over $U$, then
there is a positive integer $n$ such that
$nD$ is a Cartier divisor and
$H^0(X, \OO_X(nD)) \otimes \OO_X \to \OO_X(nD)$
is surjective on $U$.

\item
If $D$ is weakly positive over $U$, then,
for any ample $\QQ$-Cartier divisors $A$ on $X$,
there is a positive integer $n$ such that
$n(D+A)$ is a Cartier divisor and
\[
H^0(X, \OO_X(n(D+A))) \otimes \OO_X \to \OO_X(n(D+A))
\]
is surjective on $U$.
\end{enumerate}
\end{Proposition}

\section{Proof of relative Bogomolov's inequality}
Let $X$ be an algebraic variety over $k$, $x$ a point of $X$, and
$E$ a coherent $\OO_X$-module on $X$.
We say $E$ is {\em generated by global sections at $x$} if
$H^0(X, E) \otimes \OO_X \to E$ is surjective at $x$.
Let us begin with the following proposition.

\begin{Proposition}[$\ch(k) \geq 0$]
\label{prop:samp:det}
Let $X$ be a smooth algebraic variety over $k$,
$E$ a coherent $\OO_X$-module, and $x$ a point of $X$.
If $E$ is generated by global sections at $x$ and
$E$ is free at $x$, then $\det(E)$ is
generated by global sections at $x$,
where $\det(E)$ is the determinant line bundle of $E$
in the sense of \cite{KM}.
\end{Proposition}

\Proof
Let $T$ be the torsion part of $E$. Then,
$\det(E) = \det(E/T) \otimes \det(T)$.
If we set
\[
D = \sum_{\substack{P \in X, \\ \operatorname{depth}(P) = 1}} 
\operatorname{lenght}(T_P) \overline{ \{ P \} },
\]
then $\det(T) \simeq \OO_X(D)$,
where $\overline{ \{ P \} }$ is the Zariski closure of
$\{ P \}$ in $X$. Here since $E$ is free at $x$, $x \not\in \Supp(D)$.
Thus, $\det(T)$ is generated by global sections at $x$.
Moreover, it is easy to see that
$E/T$ is generated by global sections at $x$.
Therefore, to prove our proposition, we may assume that $E$ is a torsion free sheaf.

Let $r$ be the rank of $E$ and
$\kappa(x)$ the residue field of $x$.
Then, by our assumption,
there are sections $s_1, \ldots, s_r$ of $E$ such that
$\{ s_i(x) \}$ forms a basis of $E \otimes \kappa(x)$.
Thus, $s = s_1 \wedge \cdots \wedge s_r$ gives rise to
a section of $\det(E) = \left( \bigwedge^r E \right)^{**}$ with $s(y) \not= 0$.
Hence, we get our proposition.
\QED

Next let us consider the following proposition.

\begin{Proposition}[$\ch(k) \geq 0$]
\label{prop:property:global:gen}
Let $\pi : X \to Y$ be a proper and generically finite morphism of algebraic varieties over $k$.
Let $y$ be a point of $Y$ such that $\pi$ is finite over $y$.
\begin{enumerate}
\renewcommand{\labelenumi}{(\arabic{enumi})}
\item 
Let $\phi : E \to Q$ be a homomorphism of coherent $\OO_X$-modules.
If $\phi$ is surjective at each point of $\pi^{-1}(y)$
and $\pi_*(E)$ is generated by global sections
at $y$, then $\pi_*(Q)$ is generated by global sections at $y$.

\item
Let $E_1$ and $E_2$ be coherent $\OO_X$-modules.
If $\pi_*(E_1)$ and $\pi_*(E_2)$ are generated by global sections at $y$, then 
so is $\pi_*(E_1 \otimes E_2)$ at $y$.

\item
Let $E$ be a coherent $\OO_X$-module.
If $\pi_*(E)$ is generated by global sections at $y$, 
then so is $\pi_*(\Sym^n(E))$ at $y$ for every $n > 0$.
\end{enumerate}
\end{Proposition}

\Proof
(1) We can take an affine open neighborhood $U$ of $y$ such that
$\pi$ is finite over $U$ and $\phi : E \to Q$ is surjective over $\pi^{-1}(U)$.
Thus, $\pi_*(E) \to \pi_*(Q)$ is surjective at $y$.
Hence, considering the following diagram:
\[ 
\begin{CD}
H^0(Y, \pi_*(E)) \otimes \OO_Y @>>> \pi_*(E) \\
@VVV @VVV \\
H^0(Y, \pi_*(Q)) \otimes \OO_Y @>>> \pi_*(Q),
\end{CD}
\]
we have our assertion.

\medskip
(2)
Let $U$ be an affine open neighborhood of $y$ such that
$\pi$ is finite over $U$. We set $U = \Spec(A)$ for some integral domain $A$.
Since $\pi$ is finite over $U$, there is an integral domain $B$ with $\pi^{-1}(U) = \Spec(B)$.
Here we take $B$-modules $M_1$ and $M_2$ such that
$M_1$ and $M_2$ give rise to $\rest{E_1}{\pi^{-1}(U)}$ and $\rest{E_2}{\pi^{-1}(U)}$
respectively. Then, we have a natural surjective homomorphism
$M_1 \otimes_A M_2 \to M_1 \otimes_B M_2$. This shows us that
$\pi_*(E) \otimes \pi_*(E_2) \to \pi_*(E_1 \otimes E_2)$ is surjective at $y$.
Here, let us consider the following diagram:
\[ 
\begin{CD}
H^0(Y, \pi_*(E_1)) \otimes H^0(Y, \pi_*(E_2)) \otimes \OO_Y @>>> \pi_*(E_1) \otimes \pi_*(E_2) \\
@VVV @VVV \\
H^0(Y, \pi_*(E_1 \otimes E_2)) \otimes \OO_Y @>>> \pi_*(E_1 \otimes E_2),
\end{CD}
\]
where 
$H^0(Y, \pi_*(E_1)) \otimes H^0(Y, \pi_*(E_2)) \otimes \OO_Y \to \pi_*(E_1) \otimes \pi_*(E_2)$
is surjective at $y$ by our assumption.
Thus, we get (2).

\medskip
(3) This is a consequence of (1) and (2) because we have a natural surjective homomorphism
$E^{\otimes n} \to \Sym^n(E)$.
\QED

Before starting the main theorem, we need to
prepare the following formula derived from Grothendieck-Riemann-Roch theorem.

\begin{Lemma}[$\ch(k) \geq 0$]
\label{lem:growth:c1:line:by:R:R}
Let $X$ and $Y$ be algebraic varieties over $k$, and
$f : X \to Y$ a surjective and projective morphism over $k$ of $\dim f = d$.
Let $L$ and $A$ be line bundles on $X$. If $Y$ is smooth, then
there are elements $Z_1, \ldots, Z_d$ of $\Chow_{\dim Y - 1}(Y) \otimes {\QQ}$ such that
\[
c_1\left( Rf_*(L^{\otimes n} \otimes A) \right) \cap [Y] = 
\frac{ f_*(c_1(L)^{d+1} \cap [X])}{(d+1)!} n^{d+1} + 
\sum_{i=0}^{d} Z_i n^i
\]
for all $n > 0$.
\end{Lemma}

\Proof
We use the same symbol as in \cite{Fu}.
First of all,
$Rf_*(L^{\otimes n} \otimes A) \in K^{\circ}(Y)$ because $Y$ is smooth. Thus, by
\cite[Theorem~18.3, (1) and (2)]{Fu}, i.e.,
Riemann-Roch theorem for singular varieties,
\addtocounter{Claim}{1}
\begin{equation}
\label{eqn:lem:growth:c1:line:by:R:R:1}
\chern(Rf_*(L^{\otimes n} \otimes A)) \cap \tau_Y(\OO_Y) =
f_*(\chern(L^{\otimes n} \otimes A) \cap \tau_X(\OO_X)).
\end{equation}
Since $\tau_X(\OO_X) = [X] + \text{terms of dimension $< \dim X$}$ by
\cite[Theorem~18.3, (5)]{Fu}, it is easy to see that
there are $T_0, \ldots, T_{d} \in \Chow_{\dim Y - 1}(X) \otimes \QQ$ such that
\[
\left( \chern(L^{\otimes n} \otimes A) \cap \tau_X(\OO_X) \right)_{\dim Y - 1}
= \frac{c_1(L)^{d+1} \cap [X]}{(d+1)!} n^{d+1} + \sum_{i=0}^d T_i n^i.
\]
Thus,
\addtocounter{Claim}{1}
\begin{equation}
\label{eqn:lem:growth:c1:line:by:R:R:2}
f_*(\chern(L^{\otimes n} \otimes A) \cap \tau_X(\OO_X))_{\dim Y - 1} =
\frac{f_*(c_1(L)^{d+1} \cap [X])}{(d+1)!} n^{d+1} + \sum_{i=0}^d f_*(T_i) n^i.
\end{equation}
On the other hand, since
$\tau_Y(\OO_Y) = [Y] + \text{terms of dimension $< \dim Y$}$,
if we denote by $S$ the $(\dim Y - 1)$-dimensional part of $\tau_Y(\OO_Y)$,
then
\addtocounter{Claim}{1}
\begin{equation}
\label{eqn:lem:growth:c1:line:by:R:R:3}
\left( \chern(Rf_*(L^{\otimes n} \otimes A)) \cap \tau_Y(\OO_Y) \right)_{\dim Y - 1} =
c_1(Rf_*(L^{\otimes n} \otimes A)) \cap [Y] + \rank (Rf_*(L^{\otimes n} \otimes A)) S.
\end{equation}
Here, $\rank (Rf_*(L^{\otimes n} \otimes A)) = 
\chi(X_{\eta}, (L^{\otimes n} \otimes A)_{\eta})$
is a polynomial of $n$ with degree $d$ at most,
where $\eta$ is the generic point of $Y$.
Thus, combining \eqref{eqn:lem:growth:c1:line:by:R:R:1},
\eqref{eqn:lem:growth:c1:line:by:R:R:2} and \eqref{eqn:lem:growth:c1:line:by:R:R:3},
we have our lemma.
\QED

Let us start the main theorem of this paper.

\begin{Theorem}[$\ch(k) \geq 0$]
\label{thm:nef:psudo:dis}
Let $X$ be a quasi-projective variety over $k$, $Y$ a smooth quasi-projective variety over $k$,
and $f : X \to Y$ a surjective and projective morphism over $k$ of $\dim f = 1$.
Let $F$ be a locally free sheaf on $X$ with 
$f_*(c_1(F) \cap [X]) = 0$, and $y$ a point of $Y$.
We assume that
$f$ is flat over $y$, and that there are line bundles $L$ and $M$ on
the geometric fiber $X_{\bar{y}}$ over $y$
such that
\[
H^0(X_{\bar{y}}, \Sym^m(F_{\bar{y}}) \otimes L) = 
H^1(X_{\bar{y}}, \Sym^m(F_{\bar{y}}) \otimes M) = 0
\]
for $m \gg 0$.
Then, $f_*\left( (c_2(F) - c_1(F)^2) \cap [X]) \right)$ is weakly positive at $y$.
\end{Theorem}

\Proof
Let $A$ be a very ample line bundle on $X$ such that 
$A_{\bar{y}} \otimes L$ and $A_{\bar{y}} \otimes M^{\otimes -1}$ are very ample on $X_{\bar{y}}$.
First of all, we would like to see the following.

\begin{Claim}
\label{claim:thm:nef:psudo:dis:0}
$H^0(X_y, \Sym^m(F_y) \otimes A_y^{\otimes -1}) = 
H^1(X_y, \Sym^m(F_y) \otimes A_y) = 0$
for $m \gg 0$.
\end{Claim}

In general, for a coherent sheaf $G$ on $X_y$,
$H^i(X_{\bar{y}}, G \otimes_{\kappa(y)} \overline{\kappa(y)}) =
H^i(X_y, G) \otimes_{\kappa(y)} \overline{\kappa(y)}$ for all $i \geq 0$.
Thus, it is sufficient to show that
\[
H^0(X_{\bar{y}}, \Sym^m(F_{\bar{y}}) \otimes A_{\bar{y}}^{\otimes -1}) = 
H^1(X_{\bar{y}}, \Sym^m(F_{\bar{y}}) \otimes A_{\bar{y}}) = 0
\]
for $m \gg 0$.
Since $A_{\bar{y}} \otimes L$ is very ample and
$\#(\overline{\kappa(y)}) = \infty$, there is a section 
$s \in H^0(X_{\bar{y}}, A_{\bar{y}} \otimes L)$
such that $s \not= 0$ in $(A_{\bar{y}} \otimes L) \otimes \kappa(P)$
for any associated points $P$ of $X_{\bar{y}}$.
Then, $\OO_{X_{\bar{y}}} \overset{\times s}{\longrightarrow} A_{\bar{y}} \otimes L$
is injective.
Thus, tensoring the above injection with
$\Sym^m(F_{\bar{y}}) \otimes A_{\bar{y}}^{\otimes -1}$, we have an injection
\[
\Sym^m(F_{\bar{y}}) \otimes A_{\bar{y}}^{\otimes -1} \to \Sym^m(F_{\bar{y}}) \otimes L.
\]
Hence, $H^0(X_{\bar{y}}, \Sym^m(F_{\bar{y}}) \otimes A_{\bar{y}}^{\otimes -1}) = 0$ for $m \gg 0$.

In the same way,
there is a section $s' \in H^0(X_{\bar{y}}, A_{\bar{y}} \otimes M^{\otimes -1})$
such that $s' \not= 0$ in $(A_{\bar{y}} \otimes M^{\otimes -1}) \otimes \kappa(P)$
for any associated points $P$ of $X_{\bar{y}}$.
Then, $\OO_{X_{\bar{y}}} \overset{\times s'}{\longrightarrow} A_{\bar{y}} \otimes M^{\otimes -1}$
is injective and its cokernel $T$ has the $0$-dimensional support.
Thus, tensoring an exact sequence
\[
0 \to \OO_{X_{\bar{y}}} \to A_{\bar{y}} \otimes M^{\otimes -1} \to T \to 0
\]
with $\Sym^m(F_{\bar{y}}) \otimes M$, we obtain an exact sequence
\[
0 \to \Sym^m(F_{\bar{y}}) \otimes M \to \Sym^m(F_{\bar{y}}) \otimes A_{\bar{y}} \to 
\Sym^m(F_{\bar{y}}) \otimes M \otimes T \to 0.
\]
Hence, we get a surjection
\[
H^1(X_{\bar{y}}, \Sym^m(F_{\bar{y}}) \otimes M) \to H^1(X_{\bar{y}}, \Sym^m(F_{\bar{y}}) \otimes A_{\bar{y}}).
\]
Therefore, $H^1(X_{\bar{y}}, \Sym^m(F_{\bar{y}}) \otimes A_{\bar{y}}) = 0$ for $m \gg 0$.

\bigskip
Since $X$ is an integral scheme over $k$ 
of dimension greater than or equal to $2$,
and $X_y$ is a $1$-dimensional scheme over $\kappa(y)$,
by virtue of \cite[Theorem~6.10]{JB},
there is $B \in |A^{\otimes 2}|$ such that
$B$ is integral, and that $B \cap X_y$ is finite, i.e.,
$B$ is finite over $y$.
Let $\pi : B \to Y$ be the morphism induced by $f$.
Let $H$ be an ample line bundle on $Y$ such that
$\pi_*(F_B)  \otimes H$ and $\pi_*(A_B) \otimes H$ are generated by global sections at $y$,
where $F_B = \rest{F}{B}$ and $A_B = \rest{A}{B}$.

Let $\mu : P = \Proj\left( \bigoplus_{m=0}^{\infty} \Sym^m(F) \right) \to X$
be the projective bundle and
$\OO_{P}(1)$ the tautological line bundle on $P$.
We set $h = f \cdot \mu : P \to Y$.
Let us consider 
\[
c_1 \left( Rh_*((\OO_{P}(1) \otimes h^*(H))^{\otimes m} 
\otimes \mu^*(A^{\otimes -1}) \otimes h^*(H)) \right)
\cap [Y]
\]
for $m \gg 0$.  By Lemma~\ref{lem:growth:c1:line:by:R:R},
there are elements $Z_0, \ldots, Z_{r}$ of
$\Chow_{\dim Y - 1}(Y) \otimes \QQ$ such that
\begin{multline*}
c_1 \left( Rh_*((\OO_{P}(1) \otimes h^*(H))^{\otimes m} \otimes \mu^*(A^{\otimes -1}) 
\otimes h^*(H) ) \right) \cap [Y] \\
= \frac{h_*( c_1(\OO_P(1) \otimes h^*(H))^{r+1} \cap [P])}{(r+1)!} m^{r + 1} 
+ \sum_{i=0}^{r} Z_i m^i,
\end{multline*}
where $r$ is the rank of $F$.
Here
\[
\begin{cases}
\mu_*(c_1(\OO_{P}(1))^{r+1} \cap [P]) = (c_1(F)^2 - c_2(F)) \cap [X], \\
\mu_*(c_1(\OO_{P}(1))^{r} \cap [P]) = c_1(F) \cap [X], \\
\mu_*(c_1(\OO_{P}(1))^{r-1} \cap [P]) = [X], \\
\mu_*(c_1(\OO_{P}(1))^{j} \cap [P]) = 0 \qquad (0 \leq j < r-1).
\end{cases}
\]
Thus, by using projection formula, we have
\begin{align*}
h_*( c_1(\OO_P(1) \otimes h^*(H))^{r+1} \cap [P]) & =
f_* \mu_*\left( 
\sum_{i=0}^{r+1} \mu^* f^* (c_1(H)^i) \cap (c_1(\OO_P(1))^{r+1-i} \cap [P]) \right) \\
& = f_* \left( (c_1(F)^2 - c_2(F)) \cap [X] \right) \\
& \qquad\qquad + r f_* \left( f^*(c_1(H)) \cap (c_1(F) \cap [X]) \right) \\
& \qquad\qquad\qquad\quad + \frac{r(r+1)}{2} f_* \left( f^*(c_1(H)^2) \cap [X] \right) \\ 
& = - f_* \left( (c_2(F) - c_1(F)^2) \cap [X] \right)
\end{align*}
because $f_* (c_1(F) \cap [X]) = 0$ and $f_*([X]) = 0$.
Moreover, 
\[
R\mu_*((\OO_{P}(1) \otimes h^*(H))^{\otimes m} \otimes \mu^*(A^{\otimes -1}) \otimes h^*(H)) =
\Sym^m(F \otimes f^*(H)) \otimes A^{\otimes -1} \otimes f^*(H).
\]
Therefore, we get
\begin{multline*}
\sum_{i \geq 0}
(-1)^i c_1\left( R^if_*\left(\Sym^{m}(F \otimes f^*(H)) \otimes A^{\otimes -1} \otimes f^*(H) 
\right) \right) \cap [Y] \\
= -\frac{1}{(r+1)!} f_* \left( (c_2(F) - c_1(F)^2) \cap [X] \right) m^{r + 1}
+ \sum_{i=0}^{r} Z_i m^i.
\end{multline*}

Here we claim the following.

\begin{Claim}
\label{claim:thm:nef:psudo:dis:1}
If $m \gg 0$, then we have the following.
\begin{enumerate}
\renewcommand{\labelenumi}{(\alph{enumi})}
\item
$c_1\left( R^if_*\left(\Sym^{m}(F \otimes f^*(H)) \otimes A^{\otimes -1} \otimes f^*(H)
\right) \right) \cap [Y] = 0$
for all $i \geq 2$.

\item
$f_*\left(\Sym^{m}(F \otimes f^*(H)) \otimes A^{\otimes -1} \otimes f^*(H) \right) = 0$.

\item
$R^1f_*\left(\Sym^{m}(F \otimes f^*(H)) \otimes A^{\otimes -1} \otimes f^*(H) \right)$
is free at $y$.

\item
$R^1 f_*\left(\Sym^{m}(F \otimes f^*(H)) \otimes A \otimes f^*(H) \right) = 0$ around $y$.
\end{enumerate}
\end{Claim}

(a) : Let $Y'$ be the maximal open set of $Y$ such that $f$ is flat over $Y'$.
If $i \geq 2$, then the support of 
$R^if_*\left(\Sym^{m}(F \otimes f^*(H) ) \otimes A^{\otimes -1} \otimes f^*(H) \right)$
is contained in $Y \setminus Y'$. Here $\codim(Y \setminus Y') \geq 2$.
Thus, we get (a).

(b) and (c) :
By Claim~\ref{claim:thm:nef:psudo:dis:0}, 
$H^0(X_y, \Sym^m(F_{y}) \otimes A_{y}^{\otimes -1}) = 0$
for $m \gg 0$. Thus, using the upper-semicontinuity of dimension of cohomology groups,
there is an open neighborhood $U_m$ of $y$
such that $f$ is flat over $U_m$ and
$H^0(X_{y'}, \Sym^m(F_{y'}) \otimes A_{y'}^{\otimes -1}) = 0$
for all $y' \in U_m$, which implies (b) because
$f_*\left(\Sym^{m}(F \otimes f^*(H)) \otimes A^{\otimes -1} \otimes f^*(H) \right)$
is torsion free.
Here, since $f$ is flat over $U_m$, $\chi(X_{y'}, \Sym^m(F_{y'}) \otimes A_{y'}^{\otimes -1})$
is a constant with respect to $y' \in U_m$. Therefore,
so is $h^1(X_{y'}, \Sym^m(F_{y'}) \otimes A_{y'}^{\otimes -1})$
with respect to $y' \in U_m$.
Thus, we have (c).

(d) : 
By virtue of Claim~\ref{claim:thm:nef:psudo:dis:0},
$H^1(X_{y}, \Sym^m(F_{y}) \otimes A_{y}) = 0$ for $m \gg 0$.
Thus, there is an open neighborhood $U'_m$ of $y$
such that $f$ is flat over $U'_m$ and
\[
H^1(X_{y'}, \Sym^m(F_{y'}) \otimes A_{y'}) = 0
\]
for all $y' \in U'_m$.
Hence, we can see (d).

\bigskip
By (a) and (b) of Claim~\ref{claim:thm:nef:psudo:dis:1},
\begin{multline*}
\frac{1}{(r+1)!} f_* \left( (c_2(F) - c_1(F)^2) \cap [X] \right) \\
= \frac{c_1\left( R^1 f_*\left(\Sym^{m}(F \otimes f^*(H)) \otimes A^{\otimes -1} \otimes f^*(H)
\right) \right) \cap [Y]}{m^{r+1}}
+ \sum_{i=0}^{r} \frac{Z_i}{m^{r+1-i}}.
\end{multline*}
Hence, it is sufficient to show that 
\[
c_1\left( R^1 f_*\left(\Sym^{m}(F \otimes f^*(H)) \otimes A^{\otimes -1} \otimes f^*(H)
\right) \right) \cap [Y]
\]
is semi-ample at $y$.

\medskip
Since $\pi_*(F_B \otimes \pi^*(H))$ and
$\pi_*(A_B \otimes \pi^*(H))$ are generated by global sections at $y$,
by (2) and (3) of Proposition~\ref{prop:property:global:gen},
$\pi_*(\Sym^m(F_B \otimes \pi^*(H)) \otimes A_B \otimes \pi^*(H))$
is generated by global sections at $y$.
On the other hand,
a short exact sequence
\begin{multline*}
0 \to \Sym^m(F \otimes f^*(H)) \otimes A^{\otimes -1} \otimes f^*(H) \to
\Sym^m(F \otimes f^*(H)) \otimes A \otimes f^*(H) \\
\to \Sym^m(F_B \otimes \pi^*(H)) \otimes A_B \otimes \pi^*(H) \to 0
\end{multline*}
gives rise to an exact sequence
\begin{multline*}
0 \to f_*(\Sym^m(F \otimes f^*(H)) \otimes A \otimes f^*(H) ) \to 
\pi_*(\Sym^m(F_B \otimes \pi^*(H)) \otimes A_B\otimes \pi^*(H)) \\
\to R^1f_*(\Sym^m(F \otimes f^*(H)) \otimes A^{\otimes -1} \otimes f^*(H))
\to R^1f_*(\Sym^m(F \otimes f^*(H)) \otimes A \otimes f^*(H)).
\end{multline*}
Thus, by  (d) of Claim~\ref{claim:thm:nef:psudo:dis:1},
\[
\phi_m : \pi_*(\Sym^m(F_B \otimes \pi^*(H)) \otimes A_B \otimes \pi^*(H)) \to 
 R^1f_*(\Sym^m(F \otimes f^*(H)) \otimes A^{\otimes -1} \otimes f^*(H))
\]
is surjective at $y$.
Therefore, applying (1) of Proposition~\ref{prop:property:global:gen} to
the case where $\operatorname{id}_Y: Y \to Y$ and $\phi = \phi_m$,
$R^1f_*(\Sym^m(F \otimes f^*(H)) \otimes A^{\otimes -1} \otimes f^*(H))$
is generated by global sections at $y$. 
Moreover, by virtue of (c) of Claim~\ref{claim:thm:nef:psudo:dis:1},
$R^1f_*(\Sym^m(F \otimes f^*(H)) \otimes A^{\otimes -1} \otimes f^*(H))$ is free at $y$.
Hence, by Proposition~\ref{prop:samp:det},
$c_1( R^1f_*(\Sym^m(F \otimes f^*(H)) \otimes A^{\otimes -1} \otimes f^*(H))) \cap [Y]$
is semi-ample at $y$.
\QED

As a corollary of Theorem~\ref{thm:nef:psudo:dis},
we have the following.

\begin{Corollary}[$\ch(k) = 0$]
\label{cor:nef:psudo:dis}
Let $X$ be a quasi-projective variety over $k$, $Y$ a smooth quasi-projective variety over $k$,
and $f : X \to Y$ a surjective and projective morphism over $k$ with $\dim f = 1$.
Let $E$ be a locally free sheaf on $X$ and
$y$ a point of $Y$.
If $f$ is flat over $y$, 
the geometric fiber $X_{\bar{y}}$ over $y$ is reduced and Gorenstein,
and $E$ is semistable on each
connected component of the normalization of
$X_{\bar{y}}$, 
then $\dis_{X/Y}(E)$ is weakly positive at $y$.
\end{Corollary}

\Proof
We set $F = \End(E)$. First, we claim the following.

\begin{Claim}
\label{claim:cor:nef:psudo:dis}
$H^0(X_{\bar{y}}, \Sym^m(F_{\bar{y}}) \otimes A^{\otimes -1}) = 0$
for any ample line bundles $A$ on $X_{\bar{y}}$ and
any $m \geq 0$.
\end{Claim}

Let $\pi : Z \to X_{\bar{y}}$ be the normalization of $X_{\bar{y}}$.
The semistability of tensor products of semistable vector bundles
in characteristic zero was studied by a lot of authors
\cite{Gi1}, \cite{Ha1}, \cite{Mi}, \cite{Ma1} and etc.
(You can find a new elementary algebraic proof in
\S\ref{sec:relative:bogomolov:inequality:positive:characteristic}, 
which works in any characteristic under strong semistability.)
Thus, by virtue of our assumption,
$\Sym^m(\pi^*(F_{\bar{y}}))$ is semistable and
of degree $0$ on each connected component of $Z$.
Hence,
\[
H^0(Z, \pi^*( \Sym^m(F_{\bar{y}}) \otimes A^{\otimes -1} )) = 0.
\]
Here, since $\OO_{X_{\bar{y}}} \to \pi_*(\OO_Z)$ is injective,
the above implies our claim.

\bigskip
Let $L$ be an ample line bundle on $X_{\bar{y}}$ such that
$L \otimes \omega_{X_{\bar{y}}}^{\otimes -1}$ is ample.
Here, since $F_{\bar{y}}^* = F_{\bar{y}}$, 
by using Serre's duality theorem,
$H^1(X_{\bar{y}}, \Sym^m(F_{\bar{y}}) \otimes L)$
is isomorphic to the dual space of
$H^0(X_{\bar{y}}, \Sym^m(F_{\bar{y}}) \otimes 
(L \otimes \omega_{X_{\bar{y}}}^{\otimes -1})^{\otimes -1})$.
Thus, by Claim~\ref{claim:cor:nef:psudo:dis},
\[
H^0(X_{\bar{y}}, \Sym^m(F_{\bar{y}}) \otimes L^{\otimes -1}) =
H^1(X_{\bar{y}}, \Sym^m(F_{\bar{y}}) \otimes L) = 0
\]
for all $m \geq 0$.
Hence, Theorem~\ref{thm:nef:psudo:dis} implies our corollary
because $c_1(F) = 0$ and $c_2(F) = 2 \rank(E) c_2(E) - (\rank(E) -1) c_1(E)^2$.
\QED

\begin{Remark}
Even if $\rank(F) = 1$, Theorem~\ref{thm:nef:psudo:dis} is
a non-trivial fact. For, if $f : X \to Y$ is
a smooth surface fibred over a projective curve, then
the assertion of it is nothing more than
the Hodge index theorem.
\end{Remark}

\section{A weakly positive divisor on the moduli space of stable curves}

Throughout this section, we assume that $\ch(k) = 0$.

Fix an integer $g \geq 2$ and a polynomial $P_g(n) = (6n-1)(g-1)$.
Let $H_g \subset \Hilb^{P_g}_{\PP^{5g-6}}$ be a subscheme
of all tri-canonically embedded stable curves over $k$, 
$Z_g \subset H_g \times \PP^{5g-6}$ 
the universal tri-canonically embedded stable curves over $k$,
and $\pi : Z_g \to H_g$ the natural projection.
Let $\Delta$ be the minimal closed subset of $H_g$ such that
$\pi$ is not smooth over a point of $\Delta$.
Then, by \cite[Theorem~(1.6) and Corollary~(1.9)]{DM},
$Z_g$ and $H_g$ are quasi-projective and smooth over $k$, 
and $\Delta$ is a divisor with only normal crossings.
Let $\Delta = \Delta_0 \cup \cdots \cup \Delta_{[g/2]}$
be the irreducible decomposition of $\Delta$ such that,
if $x \in \Delta_i \setminus \Sing(\Delta)$, then
$\pi^{-1}(x)$ is a stable curve with one node of type $i$.
We set $U = H_g \setminus \Delta$,
$H_g^0 = H_g \setminus \Sing(\Delta_1 + \cdots +\Delta_{[g/2]})$ and
$Z_g^0 = \pi^{-1}(H_g^0)$.
In \cite[\S3]{Mo5}, we constructed a reflexive sheaf $F$ on $Z_g$ with the following
properties.
\begin{enumerate}
\renewcommand{\labelenumi}{(\arabic{enumi})}
\item
$F$ is locally free on $Z_g^0$.

\item
For each $y \in H_g \setminus (\Delta_1 \cup \cdots \cup \Delta_{[g/2]})$, 
$\rest{F}{\pi^{-1}(y)} = 
\Ker\left(H^0(\omega_{\pi^{-1}(y)}) \otimes \OO_{\pi^{-1}(y)} \to
\omega_{\pi^{-1}(y)}\right)$.

\item
$\dis_{Z_g/H_g}(F) = (8g+4) \det(\pi_*(\omega_{Z_g/H_g})) - 
g \Delta_0 - \sum_{i=1}^{\left[\frac{g}{2}\right]} 4i(g-i) \Delta_i$.
\end{enumerate}
Actually, $F$ can be constructed as follows.
First of all, we set
\[
E = \Ker\left( \pi^*(\pi_*(\omega_{Z_g/H_g})) \to 
\omega_{Z_g/H_g}\right).
\]
We would like to modify $E$ along singular fibers so that
we can get our desired $F$.
For this purpose, we consider $E^0 = \rest{E}{Z_g^0}$.
It is easy to see that $E^0$ is a locally free sheaf on $Z_g^0$.
For each $i \geq 0$, we denote $\Delta_i\cap H_g^0$ by $\Delta_i^0$.
If $i \geq 1$, then there is the irreducible decomposition
$\pi^{-1}(\Delta_i^0) = C_i^1 \cup C_i^2$ 
such that the generic fiber of $\rest{\pi}{C_i^1} : C_i^1\to \Delta_i^0$
(resp. $\rest{\pi}{C_i^2} : C_i^2 \to \Delta_i^0$) is of genus $i$
(resp. $g-i$).
Moreover, if we set
\[
  Q_i^j = \Ker\left(
  \left(\rest{\pi}{C_i^j}\right)^*\left(\rest{\pi}{C_i^j}\right)_*
  (\omega_{C_i^j/\Delta_i^0})
  \longrightarrow \omega_{C_i^j/\Delta_i^0}\right)
\]
for each $i \geq 1$ and $j = 1, 2$,
then there is a natural surjective homomorphism
\[
\alpha_i^j : \rest{E^0}{C_i^j} \to Q_i^j.
\]
Here let us consider
\[
F^0= \Ker\left( \bigoplus_{i=1}^{\left[\frac{g}{2}\right]} 
\left( \alpha_i^1 \oplus \alpha_i^2 \right)
\ : \  E^0 \longrightarrow 
\bigoplus_{i=1}^{\left[\frac{g}{2}\right]} \left( Q_i^1 \oplus Q_i^2 \right) \right).
\]
As we showed in \cite[\S3]{Mo5}, 
$F^0$ is a locally free sheaf on $Z_g^0$ with
\[
\dis_{Z_g^0/H_g^0}(F^0) = (8g+4) \det(\pi_*(\omega_{Z_g^0/H_g^0})) - 
g \Delta_0^0 - \sum_{i=1}^{\left[\frac{g}{2}\right]} 4i(g-i) \Delta_i^0.
\]
Let $\nu : Z_g^0 \to Z_g$ be the natural inclusion map. Then
$F$ can be defined by $\nu_*(F^0)$.
In order to see (2), note that
$E = F$ over $H_g \setminus (\Delta_1 \cup \cdots \cup \Delta_{[g/2]})$ and
$\pi^*(\pi_*(\omega_{Z_g/H_g})) \to 
\omega_{Z_g/H_g}$ is surjective on 
$H_g \setminus (\Delta_1 \cup \cdots \cup \Delta_{[g/2]})$
(cf. \cite[Proposition~2.1.3]{Mo5}).

\bigskip
Let $\overline{\mathcal{M}}_g$ (resp. $\mathcal{M}_g$) be the moduli space of 
stable (resp. smooth) curves of genus $g$ over $k$.
Let $\phi : H_g \to \overline{\mathcal{M}}_g$ be the canonical morphism.
Let $\lambda, \delta_0, \ldots, \delta_{[g/2]} \in 
\Pic(\overline{\mathcal{M}}_g) \otimes \QQ$ such that
$\phi^*(\lambda) = \det(\pi_*(\omega_{Z_g/H_g}))$ and
$\phi^*(\delta_i) = \Delta_i$ for all $0 \leq i \leq [g/2]$.
Let us begin with the following lemma.

\begin{Lemma}
\label{lem:criterion:pamp:on:M:g}
Let $D$ be a $\QQ$-Cartier divisor on $\overline{\mathcal{M}}_g$ and
$x$ a point of $\overline{\mathcal{M}}_g$. 
If $\phi^*(D)$ is weakly positive at any points of $\phi^{-1}(x)$, then
$D$ is weakly positive at $x$.
\end{Lemma}

\Proof
It is well known that there are a surjective finite morphism
$\pi : Y \to \overline{\mathcal{M}}_g$ of normal projective varieties
and a stable curve $f : X \to Y$ of genus $g$ such that
the induced morphism $Y \to \overline{\mathcal{M}}_g$ by $f : X \to Y$
is $\pi$. Since $\pi_*(\pi^*(D)) = \deg(\pi) D$,
by Proposition~\ref{prop:samp:pamp:push},
it is sufficient to show that $\pi^*(D)$ is weakly positive
at any points of $\pi^{-1}(x)$.

Let $y$ be a point of $\pi^{-1}(x)$. Then, there is a Zariski
open neighborhood $U$ of $y$ such that
$\rest{f_*(\omega_{X/Y}^{\otimes 3})}{U}$ is free.
Thus,
\[
\Proj\left( \bigoplus_{n=0}^{\infty} \Sym^n
\left( \rest{f_*(\omega_{X/Y}^{\otimes 3})}{U} \right) \right)
\simeq U \times \PP^{5g-6}.
\]
Therefore, there is a morphism $\mu : U \to H_g$ with $\rest{\pi}{U} = \phi \cdot \mu$.
By abuse of notation,
the induced rational map $Y \dashrightarrow H_g$ is denoted by $\mu$.
Let $\nu : Y' \to Y$ be a proper birational morphism of
normal projective varieties such that
$\mu' = \mu \cdot \nu : Y' \to H_g$ is a morphism and $\nu$ is an isomorphism over
$\nu^{-1}(U)$. Then, we have the following diagram:
\[
\begin{CD}
Y' @>{\nu}>> Y \\
@V{\mu'}VV @VV{\pi}V \\
H_g @>>{\phi}> \overline{\mathcal{M}}_g.
\end{CD}
\]
This diagram is commutative because $\phi \cdot \mu' = \pi \cdot \nu$
over $\nu^{-1}(U)$.
Hence, $\nu^*(\pi^*(D)) = {\mu'}^*(\phi^*(D))$.
Moreover, $\nu_*(\nu^*(\pi^*(D))) = \pi^*(D)$.
Thus, by virtue of Proposition~\ref{prop:samp:pamp:push},
in order to see that $\pi^*(D)$ is weakly
positive at $y$, it is sufficient to check that ${\mu'}^*(\phi^*(D))$ is weakly positive at
$y \in \nu^{-1}(U)$.

By our assumption, $\phi^*(D)$ is weakly positive at $\mu'(y)$
because $\phi(\mu'(y)) = x$.
Hence, by Proposition~\ref{prop:samp:pamp:pullback},
${\mu'}^*(\phi^*(D))$ is weakly positive at $y$.
\QED

\begin{Theorem}
\label{thm:wpos:at:M:g}
$(8g+4)\lambda - g \delta_0 - \sum_{i=1}^{[g/2]} 4i(g-i) \delta_i$
is weakly positive over $\mathcal{M}_g$.
In particular, it is pseudo-effective, and numerically effective over $\mathcal{M}_g$.
\end{Theorem}

\Proof
Let $y$ be a point of $U = H_g \setminus \Delta$. 
By virtue of \cite{PR}, $\rest{F}{\pi^{-1}(y)}$ is semistable.
Thus, by Corollary~\ref{cor:nef:psudo:dis},
$\dis_{Z_g^0/H_g^0}(F^0)$ is weakly positive at $y$.
Hence, by Proposition~\ref{prop:wp:codim:2},
so is $\dis_{Z_g/H_g}(F)$ at $y$ because $\codim(H_g \setminus H_g^0) = 2$.
Thus, $\dis_{Z_g/H_g}(F)$ is weakly positive
over $U = \phi^{-1}(\mathcal{M}_g)$.
Therefore, by virtue of Lemma~\ref{lem:criterion:pamp:on:M:g},
we get our theorem.
\QED

As a corollary, we have the following.

\begin{Corollary}
\label{cor:sharp:slope:inq}
Let $X$ be a smooth projective surface over $k$, 
$Y$ a smooth projective curve over $k$, and
$f : X \to Y$ a semistable curve of genus $g \geq 2$ over $Y$.
Then, we have the inequality
\[
(8g+4) \deg(f_*(\omega_{X/Y})) \geq
g \delta_0(X/Y) + \sum_{i=1}^{[g/2]} 4 i (g-i) \delta_i(X/Y),
\]
where $\delta_i(X/Y)$ is the number of nodes of type $i$
in all singular fibers of $f$.
\end{Corollary}

\begin{Remark}
We don't know the proof of Corollary~\ref{cor:sharp:slope:inq}
without using the moduli space $\overline{\mathcal{M}}_g$.
Let $\mu : Y \to \overline{\mathcal{M}}_g$ be the morphism
induced by $f : X \to Y$.
Then, $\mu(Y)$ might pass through $\overline{\mathcal{M}}_g \setminus
\phi(H_g^0)$.
In this case, analyses of singular fibers only in $X$
seem to be very complicated.
\end{Remark}

\section{Cones of positive divisors on the moduli space of stable curves}
\label{sec:cone:positive:divisor:moduli:spacc:stable:curve}
Throughout this section, we assume that $\ch(k) = 0$.

Let $X$ be a projective variety over $k$ and
$\mathcal{C}$ a certain family of complete irreducible curves on $X$.
A $\QQ$-Cartier divisor $D$ on $X$ is said to be {\em numerically effective for $\mathcal{C}$}
if $(D \cdot C) \geq 0$ for all $C \in \mathcal{C}$.
We set
\[
\Nef(X, \mathcal{C}) = \left\{ D \in \NS(X) \otimes \QQ \mid 
\text{$D$ is numerically effective for $\mathcal{C}$} \right\}.
\]
Moreover, for subsets $A$ and $B$ in $X$, we denote by $\Curve^A_B$
the set of all irreducible complete curves $C$ on $X$
with $C \subseteq A$ and $C \cap B \not= \emptyset$.

Let $g$ be an integer greater than or equal to $2$,
$\mathcal{I}_g$ the locus of hyperelliptic curves in $\mathcal{M}_g$,
$\overline{\mathcal{I}}_g$ the closure in $\overline{\mathcal{M}}_g$, and
$\overline{\mathcal{M}}_g^{one}$ the set of all stable curves with at most one node,
i.e., if we use the notation in the previous section,
\[
\overline{\mathcal{M}}_g^{one} = \phi\left( H_g \setminus \Sing(\Delta_0 + \cdots + \Delta_{[g/2]}) \right).
\]
Let us begin with the following lemma.

\begin{Lemma}
\label{lem:existence:curve}
There are complete irreducible curves $C, C_0, \ldots, C_{[g/2]}$ on $\overline{\mathcal{M}}_g$
with the following properties.
\begin{enumerate}
\renewcommand{\labelenumi}{(\arabic{enumi})}
\item
$C, C_0, \ldots, C_{[g/2]} \in \Curve^{\overline{\mathcal{M}}_g^{one}}_{\mathcal{I}_g}$.

\item
$C \subset \mathcal{M}_g$.

\item
$C_i \subset \overline{\mathcal{I}}_g$ for all $0 \leq i \leq [g/2]$.

\item
For all $0 \leq i, j \leq [g/2]$, $(\delta_i \cdot C_j)$ is positive if $i=j$, and
$(\delta_i \cdot C_j) = 0$ if $i \not= j$.
\end{enumerate}
\end{Lemma}

\Proof
Let $\overline{\mathcal{M}}^s_g$ be Satake's compactification of
$\mathcal{M}_g$. Then, $\overline{\mathcal{M}}^s_g$ is projective and
$\codim(\overline{\mathcal{M}}^s_g \setminus \mathcal{M}_g) \geq 2$.
Pick up one point $P \in \mathcal{I}_g$. If we take general hyperplane sections
$H_1, \ldots, H_{3g-4}$ passing through $P$, then $C = H_1 \cap \ldots \cap H_{3g-4}$
is a complete irreducible curve with $C \subseteq \mathcal{M}_g$ and $P \in C$.

Applying Proposition~\ref{prop:hyperelliptic:fibration:2} to the case 
where $a=0$, and
contracting all $(-2)$-curves in all singular fibers, we have
a stable fibred surface $f_0 : X_0 \to Y_0$ such that
$Y_0$ is projective, the generic fiber of $f_0$ is a smooth hyperelliptic curve of genus $g$,
$f_0$ is not smooth, and that
every singular fiber of $f_0$ is an irreducible nodal curve with one node.
Let $\mu_0 : Y_0 \to \overline{\mathcal{M}}_g$ be the induced morphism by $f_0 : X_0 \to Y_0$.
Then, $C_0 = \mu(Y_0)$ is our desired curve.

Finally, we fix $i$ with $1 \leq i \leq [g/2]$.
Using Proposition~\ref{prop:hyperelliptic:fibration},
there is a stable fibred surface $f_i : X_i\to Y_i$ such that
$Y_i$ is projective, the generic fiber of $f_i$ is a smooth hyperelliptic curve of genus $g$,
$f_i$ is not smooth, and that
every singular fiber of $f_i$ is a reducible curve with one node of type $i$.
Let $\mu_i : Y_i \to \overline{\mathcal{M}}_g$ be the induced morphism by 
$f_i : X_i \to Y_i$.
If we set $C_i = \mu_i(Y_i)$, then $C_i$ satisfies our requirements.
\QED

By using curves in Lemma~\ref{lem:existence:curve}, we can show the following
proposition.

\begin{Proposition}
\label{prop:cone:nef:include:special:cone}
\[
\Nef\left(\overline{\mathcal{M}}_g, \Curve^{\overline{\mathcal{M}}_g^{one}}_{\mathcal{I}_g}\right)
\subseteq \left\{ x \lambda + \sum_{i=0}^{[g/2]} y_i \delta_i \ \left| \ 
\begin{array}{l}
x \geq 0, \\
g x + (8g + 4) y_0 \geq 0, \\
i(g-i) x + (2g+1) y_i \geq 0 \quad (1 \leq i \leq [g/2]).
\end{array}
\right. \right\}.
\]
\end{Proposition}

\Proof
Let $D = x \lambda + \sum_{i=0}^{[g/2]} y_i \delta_i$ be an arbitrary element of
$\Nef(\overline{\mathcal{M}}_g, \Curve^{\overline{\mathcal{M}}_g^{one}}_{\mathcal{I}_g})$.
Let $C, C_0, \ldots, C_{[g/2]}$ be irreducible complete curves as in Lemma~\ref{lem:existence:curve}.
Then, $0 \leq (D \cdot C) = x (\lambda \cdot C)$. Hence $x \geq 0$.

To get other inequalities, we need some facts of hyperelliptic fibrations.
Details can be found in \cite[\S4, b]{CH}.
For $i > 0$, $\Delta_i \cap \overline{\mathcal{I}}_g$ is irreducible.
$\Delta_0 \cap \overline{\mathcal{I}}_g$ is however reducible and has
$1 + [(g-1)/2]$ irreducible components, say,
$\Sigma_0, \Sigma_1, \ldots, \Sigma_{[(g-1)/2]}$.
Here a general point of $\Sigma_0$ represents an irreducible curve of one node,
and a general point of $\Sigma_i$ ($i > 0$) represents
a stable curve consisting of a curve of genus $i$ and
a curve of genus $g - i - 1$ joined at two points.
The class of $\Sigma_i$ in $\Pic(\overline{\mathcal{I}}_g) \otimes \QQ$
is denoted by $\sigma_i$, and by abuse of notation,
$\rest{\delta_i}{\overline{\mathcal{I}}_g}$ is denoted by $\delta_i$.
Further, $\rest{\lambda}{\overline{\mathcal{I}}_g}$ is denoted by $\lambda$.
Then, by virtue of \cite[Proposition~(4.7)]{CH},
\[
   \delta_0 = \sigma_0 + 2 \left(\sigma_1 + \cdots + \sigma_{[(g-1)/2]} \right)
\]
and
\[
(8g+4) \lambda = g \sigma_0 + \sum_{j=1}^{[(g-1)/2]} 2(j+1)(g-j) \sigma_j
+ \sum_{i=1}^{[g/2]} 4 i (g-i) \delta_i.
\]

Let us consider $D$ as a divisor on $\overline{\mathcal{I}}_g$. Using the above relations
between $\lambda$, $\delta_i$'s and $\sigma_j$'s,
we have
\[
D = 
\left( \frac{g}{8g + 4} x + y_0 \right) \sigma_0 +
2 \sum_{j=1}^{[(g-1)/2]} \left( \frac{(j+1)(g-j)}{8g+4} x + y_0 \right) \sigma_j +
\sum_{i=1}^{[g/2]} \left( \frac{i(g-i)}{2g+1} x + y_i \right)\delta_i.
\]
Note that $C_i \cap \Sigma_j = \emptyset$ for all $0 \leq i \leq [g/2]$
and $1 \leq j \leq [(g-1)/2]$ 
because $C_i \subset \overline{\mathcal{M}}_g^{one}$.
Thus, considering $(D \cdot C_i)$, we have the remaining inequalities.
\QED

\begin{Corollary}
\label{cor:nef:cone:equal}
If $\mathcal{C}$ is a set of 
complete irreducible curves on $\overline{\mathcal{M}}_g$ with
$\Curve^{\overline{\mathcal{M}}_g^{one}}_{\mathcal{I}_g} \subseteq \mathcal{C} \subseteq
\Curve^{\overline{\mathcal{M}}_g}_{\mathcal{M}_g}$,
then
\[
\Nef\left(\overline{\mathcal{M}}_g, \mathcal{C} \right)
= \left\{ x \lambda + \sum_{i=0}^{[g/2]} y_i \delta_i\ \left| \ 
\begin{array}{l}
x \geq 0, \\
g x + (8g + 4) y_0 \geq 0, \\
i(g-i) x + (2g+1) y_i \geq 0 \quad (1 \leq i \leq [g/2]).
\end{array}
\right. \right\}.
\]
\end{Corollary}

\Proof
Since $\Nef\left(\overline{\mathcal{M}}_g, \mathcal{C} \right) \subseteq
\Nef\left(\overline{\mathcal{M}}_g, \Curve^{\overline{\mathcal{M}}_g^{one}}_{\mathcal{I}_g}\right)$,
the direction ``$\subseteq$'' is a consequence of
Proposition~\ref{prop:cone:nef:include:special:cone}.
Conversely, we assume that $D = x \lambda + \sum_{i=0}^{[g/2]} y_i \delta_i$ satisfies
\[
\begin{cases}
x \geq 0, \\
g x + (8g + 4) y_0 \geq 0, \\
i(g-i) x + (2g+1) y_i \geq 0 \quad (1 \leq i \leq [g/2]).
\end{cases}
\]
Then, since
\begin{multline*}
D = \frac{x}{8g+4} \left(
(8g+4)\lambda - g \delta_0 - \sum_{i=1}^{[g/2]} 4i(g-i) \delta_i \right) \\
+ \left(y_0 + \frac{g}{8g+4} x\right)\delta_0 + 
\sum_{i=1}^{[g/2]} \left(y_i + \frac{i(g-i)}{2g+1} x \right) \delta_i
\end{multline*}
and $\mathcal{C} \subseteq
\Curve^{\overline{\mathcal{M}}_g}_{\mathcal{M}_g}$,
we can see that $D$ is numerically effective for $\mathcal{C}$
by using Theorem~\ref{thm:wpos:at:M:g}.
\QED

In the same way, we can see the following.

\begin{Corollary}
\label{cor:wp:cone:equal}
If we set
\[
\WP(\overline{\mathcal{M}}_g; \mathcal{M}_g) = 
\{ D \in \Pic(\overline{\mathcal{M}}_g) \otimes \QQ \mid
\text{$D$ is weakly positive over $\mathcal{M}_g$} \},
\]
then
\[
\WP(\overline{\mathcal{M}}_g; \mathcal{M}_g)
= \left\{ x \lambda + \sum_{i=0}^{[g/2]} y_i \delta_i \ \left| \ 
\begin{array}{l}
x \geq 0, \\
g x + (8g + 4) y_0 \geq 0, \\
i(g-i) x + (2g+1) y_i \geq 0 \quad (1 \leq i \leq [g/2]).
\end{array}
\right. \right\}.
\]
\end{Corollary}

\Proof
Note that
\[
\WP(\overline{\mathcal{M}}_g; \mathcal{M}_g)
\subseteq \Nef\left(\overline{\mathcal{M}}_g, 
\Curve^{\overline{\mathcal{M}}_g}_{\mathcal{M}_g}\right)
\]
and that
$(8g+4)\lambda - g \delta_0 - \sum_{i=1}^{[g/2]} 4i(g-i)$ and
$\delta_i$'s are weakly positive over $\mathcal{M}_g$.
\QED

\begin{Remark}
In general, over an open set,
weak positivity is stronger 
than numerical effectivity.
Corollary~\ref{cor:nef:cone:equal} and 
Corollary~\ref{cor:wp:cone:equal} however say us that,
on the moduli space of stable curves $\overline{\mathcal{M}}_g$,
weak positivity over $\mathcal{M}_g$ coincides
with numerical effectivity over $\mathcal{M}_g$.
\end{Remark}

\section{Effective Bogomolov's conjecture over function fields}
\label{sec:bogo:conj}

Let $X$ be a smooth projective surface over $k$, 
$Y$ a smooth projective curve over $k$,
and $f : X \to Y$ a generically smooth semistable curve 
of genus $g \geq 2$ over $Y$.
Let $K$ be the function field of $Y$, $\overline{K}$
the algebraic closure of $K$, and $C$ the generic fiber of $f$.
Let $j : C(\overline{K}) \to \Jac(C)(\overline{K})$ be the map
given by $j(x) = (2g-2)x - \omega_C$ and $\Vert\ \Vert_{NT}$
the semi-norm arising from the Neron-Tate height pairing on 
$\Jac(C)(\overline{K})$.
We set 
\[
B_C(P;r) = \left\{ x \in C(\overline{K}) \mid
\Vert j(x) - P \Vert_{NT} \leq r \right\}
\]
for $P \in \Jac(C)(\overline{K})$ and $r \geq 0$, and
\[
r_C(P) = 
\begin{cases}
-\infty & \mbox{if $\#\left(B_C(P;0)\right) = \infty$}, \\
& \\
\sup \left\{ r \geq 0 \mid \#\left(B_C(P;r)\right) < \infty \right\} &
\mbox{otherwise}.
\end{cases}
\]
An effective version of Bogomolov's conjecture claims the following.

\begin{Conjecture}[Effective Bogomolov's conjecture]
\label{conj:effective:bogomolov}
If $f$ is non-isotrivial, then 
there is an effectively calculated positive number
$r_0$ with
\[
\inf_{P \in \Jac(C)(\overline{K})} r_C(P) \geq r_0.
\]
\end{Conjecture}

Recently, Ullmo \cite{Ul} proved that
$r_C(P) > 0$ for all $P \in \Jac(C)(\overline{K})$
for the case where $K$ is a number field.
As far as we know, the problem to find an effectively calculated
$r_0$ is still open.
The meaning of ``effectively calculated'' is that
a concrete algorithm or formula to find $r_0$ is required.

Here we need a rather technical condition coming from
calculations of green functions along singular fibers.
Let $\bar{f} : \overline{X} \to Y$ be the stable model of
$f : X \to Y$. Let $X_y$ (resp. $\overline{X}_y$) 
be the singular fiber of $f$ (resp. $\bar{f}$) over $y \in Y$, and
$S_y$ the set of nodes $P$ on $\overline{X}_y$ such that
$P$ is not an intersection of two irreducible components of $\overline{X}_y$,
i.e., a singularity of an irreducible component.
Let $\pi : Z_y \to \overline{X}_y$ be the partial normalization of 
$\overline{X}_y$ at each node in $S_y$.
We say $X_y$ is a {\em tree of stable components} if
the dual graph of $Z_y$  is a tree graph.
In other words, every node of type $0$ on $\overline{X}_y$ is a singularity
of an irreducible component of $\overline{X}_y$.

As an application of Corollary~\ref{cor:sharp:slope:inq},
we get the following solution of the above conjecture,
which is a generalization of \cite[Theorem~5.2]{Mo5}.

\begin{Theorem}[$\ch(k) = 0$]
\label{thm:bogomolov:function:field}
If $f$ is not smooth and every singular fiber of $f$
is a tree of stable components, then
\[
\inf_{P \in \Jac(C)(\overline{K})} r_C(P) \geq 
\sqrt{\frac{(g-1)^2}{g(2g+1)}\left(
\frac{g-1}{3}\delta_0(X/Y) + \sum_{i=1}^{\left[\frac{g}{2}\right]}
4i(g-i)\delta_i(X/Y) \right)}.
\]
\end{Theorem}

Before starting the proof of Theorem~\ref{thm:bogomolov:function:field},
let us recall several facts of green functions on a metrized graph.
For details of metrized graphs, see \cite{Zh}.

Let $G$ be a connected metrized graph and
$D$ an $\RR$-divisor on $G$.
If $\deg(D) \not= -2$, then
there are a unique measure $\mu_{(G,D)}$ on $G$ and
a unique function $g_{(G,D)}$ on $G \times G$ 
with the following properties.
\begin{enumerate}
\renewcommand{\labelenumi}{(\alph{enumi})}
\item
${\displaystyle \int_{G} \mu_{(G,D)} = 1}$.

\item
$g_{(G,D)}(x, y)$ is symmetric and continuous on $G \times G$.

\item
For a fixed $x \in G$, $\Delta_y(g_{(G,D)}(x, y)) = \delta_x - \mu_{(G,D)}$.

\item
For a fixed $x \in G$, ${\displaystyle 
\int_G g_{(G,D)}(x, y) \mu_{(G,D)}(y) = 0}$.

\item
$g_{(G,D)}(D, y) + g_{(G,D)}(y, y)$ is a constant for all $y \in G$.
\end{enumerate}
The constant $g_{(G,D)}(D, y) + g_{(G,D)}(y, y)$ is denoted by $c(G, D)$.
Further we set
\[
  \epsilon(G, D) = 2\deg(D)c(G, D) - g_{(G,D)}(D, D).
\]
We would like to calculate the invariant
$\epsilon(G,D)$ for the metrized graph $G$ with the polarization $D$.
First of all, let us see two examples, which will be elementary pieces
for our calculations.

\begin{Example}[{cf. \cite[Lemma~3.2]{Mo3}}]
\label{exam:e:for:circle}
Let $C$ be a circle of length $l$ and $O$ a point on $C$.
Then,
\[
g_{(C,0)}(O,O) = \frac{l}{12}
\quad\text{and}\quad
\epsilon(C, 0) = 0.
\]
\end{Example}

\begin{Example}[{cf. \cite[Lemma~4.4]{Mo5}}]
\label{exam:e:for:1:segment}
Let $G$ be a segment of length $l$, and $P$ and $Q$ 
terminal points of $G$. Let $a$ and $b$ be real numbers
with $a + b \not= 0$, and $D$ an $\RR$-divisor on $G$ given by
$D = (2a-1)P + (2b-1)Q$. Then,
\[
  \epsilon(G, D) = \left(\frac{4ab}{a+b} - 1\right)l,\quad
  g_{(G,D)}(P,P) = \frac{b^2}{(a+b)^2}l \quad\text{and}\quad
  g_{(G,D)}(Q,Q) = \frac{a^2}{(a+b)^2}l.
\]
\end{Example}

Let $G_1$ and $G_2$ be metrized graphs.
Fix points $x_1 \in G_1$ and $x_2 \in G_2$.
The one point sum $G_1 \vee G_2$ with respect to $x_1$ and $x_2$,
defined by $G_1 \times \{x_2\} \cup \{x_1\} \times G_2$ in $G_1 \times G_2$,
is a metrized graph obtained by joining $x_1\in G_1$ and $x_2 \in G_2$.
The joining point, which is $\{x_1\}\times\{x_2\}$ in $G_1 \times G_2$,
is denoted by $j(G_1 \vee G_2)$.
Any $\RR$-divisor on $G_i$ ($i=1,2$) can be viewed as an $\RR$-divisor on
$G_1 \vee G_2$. Then, our basic tool for our calculations is the following.

\begin{Proposition}[{cf. \cite[Proposition~4.2]{Mo5}}]
\label{prop:e:for:join:graph}
Let $G_1$ and $G_2$ be connected metrized graphs,
and $D_1$ and $D_2$ $\RR$-divisors on $G_1$ and $G_2$ respectively
with $\deg(D_i) \not= -2$ \textup{(}$i=1,2$\textup{)}. 
Let $G = G_1 \vee G_2$, $O = j(G_1\vee G_2)$, and
$D = D_1 + D_2$ on $G_1 \vee G_2$. If $\deg(D_1 + D_2) \not= -2$, then
we have the following formulae, where $d_i = \deg(D_i)$ \rom{(}$i=1, 2$\rom{)}
and $r_{G_2}(O,P)$ is the resistance between $O$ and $P$ on $G_2$.
\begin{enumerate}
\renewcommand{\labelenumi}{(\arabic{enumi})}
\item
If $P \in G_2$, then
{\allowdisplaybreaks
\begin{align*}
g_{(G, D)}(P,P) & = 
\frac{d_1}{d_1+d_2+2} r_{G_2}(O,P) 
+\frac{d_2 + 2}{d_1+d_2+2} g_{(G_2, D_2)}(P,P) \\
& \quad -\frac{d_1(d_2+2)}
{(d_1+d_2+2)^2}g_{(G_2,D_2)}(O,O) 
+\frac{(d_1 +2)^2}{(d_1 + d_2 + 2)^2}g_{(G_1, D_1)}(O,O).
\end{align*}}

\item
\iftwelvept
\[
\epsilon(G, D) = \epsilon(G_1, D_1) + \epsilon(G_2, D_2)
+ \frac{2 d_2(d_1 + 2)g_{(G_1, D_1)}(O,O) +
2 d_1(d_2 + 2)g_{(G_2, D_2)}(O, O)}{d_1 + d_2 + 2}.
\]
\else
\begin{multline*}
\epsilon(G, D) = \epsilon(G_1, D_1) + \epsilon(G_2, D_2) \\
+ \frac{2 d_2(d_1 + 2)g_{(G_1, D_1)}(O,O) +
2 d_1(d_2 + 2)g_{(G_2, D_2)}(O, O)}{d_1 + d_2 + 2}.
\end{multline*}
\fi
\end{enumerate}
\end{Proposition}

Combining the above proposition and Example~\ref{exam:e:for:circle},
we have the following.

\begin{Corollary}
\label{cor:e:for:join:graph:circle}
Let $G$ be a connected metrized graph and $D$ an $\RR$-divisor on $G$ with
$\deg(D) \not= -2$. Let $C$ be a circle of length $l$.
Then,
\[
    \epsilon(G \vee C, D) = \epsilon(G, D) + \frac{\deg(D)}{3(\deg D + 2)} l.
\]
\end{Corollary}

Let $G$ be a connected metrized graph. We assume that $G$ is a tree, i.e.,
there is no loop in $G$.
Let $\Vt(G)$ (resp. $\Ed(G)$) be the set of vetexes (resp. edges)
of $G$.
For a function $\alpha : \Vt(G) \to \RR$, we define the divisor
$D(\alpha)$ on $G$ to be
\[
D(\alpha) = \sum_{x \in \Vt(G)} (2\alpha(x) - 2 + v(x)) x,
\]
where $v(x)$ is the number of branches starting from $x$.
It is easy to see that
\[
\deg(D(\alpha)) + 2 = 2 \sum_{x \in \Vt(G)} \alpha(x).
\]
To give an exact formula for $\epsilon(G, D(\alpha))$,
we need to introduce the following notation.
Let $e$ be an edge of $G$, $P$ and $Q$ terminal points of $e$, and
$e^{\circ} = e \setminus \{ P, Q \}$. 
Since $G$ is a connected tree, there are two connected
sub-graphs $G'_e$ and $G''_e$ such that
$G \setminus e^{\circ} = G'_e \coprod G''_e$.
Then, we have the following.

\begin{Proposition}
\label{prop:cal:e:G:D:tree}
With the same notation as above,
if $\alpha(x) \geq 0$ for all $x \in \Vt(G)$
and $\sum_{x \in \Vt(G)} \alpha(x) \not= 0$, then
\[
\epsilon(G, D(\alpha)) =
\sum_{e \in \Ed(G)} \left(
\frac{4 \left( \sum_{x \in \Vt(G'_e)} \alpha(x) \right)
\left( \sum_{x \in \Vt(G''_e)} \alpha(x) \right)}
{\sum_{x \in \Vt(G)} \alpha(x) } - 1
\right) l(e),
\]
where $l(e)$ is the length of $e$.
\end{Proposition}

\Proof
For a positive number $t$, we set
$\alpha_t(x) = \alpha(x) + t$. Then, it is easy to see that
\[
\lim_{t \downarrow 0} \epsilon(G, D(\alpha_t)) = \epsilon(G, D(\alpha)).
\]
Thus, in order to prove our proposition, we may assume that
$\alpha(x) > 0$ for all $x \in \Vt(G)$.

We fix $P \in \Vt(G)$. For $e \in \Ed(G)$,
we denote by $G_{P,e}$ the connected component of $G \setminus e^{\circ}$
not containing $P$, i.e., if $P \not\in G'_e$, then $G_{P,e} = G'_e$;
otherwise, $G_{P,e} = G''_e$.
With this notation, let us consider the following claim.

\begin{Claim}
\label{claim:cal:g:G:D:tree}
\[
g_{(G, D(\alpha))}(P, P) =
\sum_{e \in \Ed(G)} \frac{\left(\sum_{x \in \Vt(G_{P, e})} \alpha(x)\right)^2}
{\left(\sum_{x \in \Vt(G)} \alpha(x)\right)^2} l(e).
\]
\end{Claim}

We prove this claim by induction on $\#(\Ed(G))$.
If $\#(\Ed(G)) = 0, 1$, then our assertion is obvious by
Example~\ref{exam:e:for:1:segment}.
Thus, we may assume that  $\#(\Ed(G)) \geq 2$.

First, we suppose that $P$ is not a terminal point.
Let $G'$ be one branch starting from $P$, and $G''$ a connected
sub-graph such that $G' \cup G'' = G$ and $G' \cap G'' = \{ P \}$.
We define $\alpha' : \Vt(G') \to \RR$ and
$\alpha'' : \Vt(G'') \to \RR$ by
\[
\alpha'(x) =
\begin{cases}
1 & \text{if $x = P$} \\
\alpha(x) & \text{otherwise}
\end{cases}
\]
and $\alpha'' = \rest{\alpha}{\Vt(G'')}$.
Then, we have $G = G' \vee G''$ and
$D(\alpha) = D(\alpha') + D(\alpha'')$.
Thus, using (1) of Proposition~\ref{prop:e:for:join:graph} and 
hypothesis of induction,
we can easily see our claim.

Next we suppose that $P$ is a terminal point.
Pick up $e \in \Ed(G)$ such that $P$ is a terminal of $e$.
Let $O$ be another terminal of $e$.
We set $G' = e$ and $G'' = (G \setminus e) \cup \{ O \}$.
Moreover, we define $\alpha' : \Vt(G') = \{P, O\} \to \RR$ and
$\alpha'' : \Vt(G'') \to \RR$ by
$\alpha'(P) = \alpha(P)$, $\alpha'(O) = 1$ and
$\alpha'' = \rest{\alpha}{\Vt(G'')}$.
Then, $G = G' \vee G''$ and $D(\alpha) = D(\alpha') + D(\alpha'')$.
Thus, using (1) of Proposition~\ref{prop:e:for:join:graph},
Example~\ref{exam:e:for:1:segment} and
hypothesis of induction,
we can see our claim after easy calculations.

\medskip
Let us go back to the proof of Proposition~\ref{prop:cal:e:G:D:tree}.
We prove it by induction on $\#(\Ed(G))$.
If $\#(\Ed(G)) = 0, 1$, then our assertion comes from
Example~\ref{exam:e:for:1:segment}.
Thus, we may assume that $\#(\Ed(G)) \geq 2$.
Let us pick up a terminal edge $e$ of $G$.
Let $\{ O, P \}$ be terminals of $e$ such that
$P$ gives a terminal of $G$. We set $G_1 = e$ and 
$G_2 = (G \setminus e) \cup \{ O \}$.
Moreover, we define $\alpha_1 : \{ O, P \} = \Vt(G_1) \to \RR$ and
$\alpha_2 : \Vt(G_2) \to \RR$ by $\alpha_1(O) = 1$, $\alpha_1(P) = \alpha(P)$, and
$\alpha_2 = \rest{\alpha}{\Vt(G_2)}$.
Then, $G = G_1 \vee G_2$ and $D(\alpha) = D(\alpha_1) + D(\alpha_2)$.
Thus, if we set
\[
\begin{cases}
A = \sum_{x \in \Vt(G)} \alpha(x), \\
a = \alpha(P), \\
A_{e'} = \sum_{x \in \Vt(G_{O, e'})} \alpha(x) & 
\text{for $e' \in \Vt(G) \setminus \{ e \}$},
\end{cases}
\]
then, by (2) of Proposition~\ref{prop:e:for:join:graph},
Example~\ref{exam:e:for:1:segment}, Claim~\ref{claim:cal:g:G:D:tree}
and hypothesis of induction,
we have
{\allowdisplaybreaks
\begin{align*}
\epsilon(G, D(\alpha)) & =
\left(\frac{4a}{a+1} - 1 \right) l(e) + 
\sum_{e' \in \Vt(G) \setminus \{ e \}}
\left(\frac{4 A_{e'}(A - a - A_{e'})}{A - a} - 1 \right) l(e') \\
& \qquad 
+ \frac{4(A - a - 1)(a+1)}{A} \frac{a^2}{(a+1)^2} l(e)
+ \frac{4a(A-a)}{A} \sum_{e' \in \Vt(G) \setminus \{ e \}}
\frac{A_{e'}^2}{(A - a)^2} l(e') \\
& = \left( \frac{4a(A-a)}{A} - 1 \right) l(e) +
\sum_{e' \in \Vt(G) \setminus \{ e \}}
\left( \frac{4 A_{e'} (A - A_{e'})}{A} - 1 \right) l(e').
\end{align*}}
Therefore, we get our proposition.
\QED

\begin{Corollary}[$\ch(k) \geq 0$]
\label{cor:e:for:semistable:chain}
Let $X$ be a smooth projective surface over $k$, 
$Y$ a smooth projective curve over $k$, and
$f : X \to Y$ a generically smooth semistable curve of genus $g \geq 2$
over $Y$.
Let $X_y$ be the singular fiber of $f$ over $y \in Y$, and
$X_y = C_1 + \cdots + C_n$ the irreducible decomposition of $X_y$.
Let $G_y$ be the metrized graph given by the configuration of $X_y$,
$v_i$ the vertex of $G_y$ corresponding to $C_i$, and
$\omega_y$ the divisor on $G_y$ defined by
$\omega_y = \sum_i (\omega_{X/Y} \cdot C_i) v_i$.
If $X_y$ is a tree of stable components, then
\[
\epsilon(G_y, \omega_y) = \frac{g-1}{3g} \delta_{0}(X_y) +
\sum_{i=1}^{\left[\frac{g}{2}\right]} 
\left( \frac{4 i (g-i)}{g} - 1 \right) \delta_{i}(X_y).
\]
\end{Corollary}

\Proof
Let $\bar{f} : \overline{X} \to Y$ be the stable model of
$f : X \to Y$, and
$S_y$ the set of nodes $P$ on $\overline{X}_y$ such that
$P$ is a singularity of an irreducible component.
Let $\pi : Z_y \to \overline{X}_y$ be the partial normalization of 
$\overline{X}_y$ at each node in $S_y$.
Let $\overline{G}_y$ be the dual graph of $Z_y$.
Let $l_1, \cdots, l_{r}$ be circles in $G_y$
corresponding to nodes in $S_y$.
Then, $G_y = \overline{G}_y \vee l_1 \vee \cdots \vee l_{r}$.
Moreover, if $g_i$ is the arithmetic genus of $C_i$ and
$\alpha : \Vt(G_y) \to \RR$ is given by
$\alpha(v_i) = g_i$, then $\omega_y = D(\alpha)$.
Here, by virtue of Proposition~\ref{prop:cal:e:G:D:tree},
\[
\epsilon(\overline{G}_y, \omega_y) =
\sum_{i=1}^{\left[\frac{g}{2}\right]} 
\left( \frac{4 i (g-i)}{g} - 1 \right) \delta_{i}(X_y).
\]
Therefore, it follows from Corollary~\ref{cor:e:for:join:graph:circle} that
\[
\epsilon(G_y, \omega_y) = \frac{g-1}{3g} \delta_{0}(X_y) +
\sum_{i=1}^{\left[\frac{g}{2}\right]} 
\left( \frac{4 i (g-i)}{g} - 1 \right) \delta_{i}(X_y).
\]
\QED

Let us start the proof of Theorem~\ref{thm:bogomolov:function:field}.
First of all, note the following fact
(cf. \cite[Theorem 5.6]{Zh}, 
\cite[Corollary 2.3]{Mo3} or \cite[Theorem 2.1]{Mo4}).
If $(\omega_{X/Y}^a \cdot \omega_{X/Y}^a)_a > 0$, then
\[
\inf_{P \in \Jac(C)(\overline{K})} r_C(P) \geq 
\sqrt{(g-1)(\omega_{X/Y}^a \cdot \omega_{X/Y}^a)_a},
\]
where $(\ \cdot \ )_a$ is the admissible pairing.

By the definition of admissible pairing,
we can see
\[
\left(\omega_{X/Y}^a \cdot \omega_{X/Y}^a \right)_a =
\left(\omega_{X/Y} \cdot \omega_{X/Y} \right) - 
\sum_{y \in Y} \epsilon(G_y, \omega_y).
\]
On the other hand,
by Corollary~\ref{cor:sharp:slope:inq},
we have
\[
(8g + 4) \deg(f_*(\omega_{X/Y})) \geq
g \delta_0(X/Y) + \sum_{i=1}^{\left[\frac{g}{2}\right]} 4i(g-i) \delta_i(X/Y).
\]
Thus, using Noether formula, the above inequality implies
\[
(\omega_{X/Y} \cdot \omega_{X/Y}) \geq
\frac{g-1}{2g+1} \delta_0(X/Y) + \sum_{i=1}^{\left[\frac{g}{2}\right]}
\left(\frac{12i(g-i)}{2g+1} - 1 \right) \delta_i(X/Y).
\]
Therefore,
we have our theorem by Corollary~\ref{cor:e:for:semistable:chain}.
\QED

Moreover, using Ullmo's result \cite{Ul} and
Proposition~\ref{prop:cal:e:G:D:tree}, we have the following.

\begin{Corollary}
Let $K$ be a number field, $O_{K}$ the ring of integers, and
$f : X \to \Spec(O_{K})$ a regular semistable arithmetic surface
of genus $g \geq 2$ over $O_{K}$. Let $S$ be the subset of
$\Spec(O_{K}) \setminus \{ 0 \}$ such that $P \in S$ if and only if
the stable model of the geometric fiber $X_{\bar{P}}$ at $P$
is a tree of stable components. Then, we have
\[
\left( \omega^{Ar}_{X/O_{K}} \cdot \omega^{Ar}_{X/O_{K}} \right) >
\sum_{P \in S} \left\{
\frac{g-1}{3g}\delta_{0}(X_{\bar{P}}) +
\sum_{i=1}^{[g/2]} \left(\frac{4i(g-i)}{g} - 1 \right)\delta_{i}(X_{\bar{P}})
\right\} \log\#(O_{K}/P).
\]
\end{Corollary}

\section{Generalization to higher dimensional fibrations}
In this section,
we consider a generalization of Corollary~\ref{cor:nef:psudo:dis}
to  higher dimensional fibrations.

First of all, let us recall the definition of semistability of vector bundles.
Let $V$ be a smooth projective variety of dimension $d$ over $k$,
and $H_1, \ldots, H_{d-1}$ ample line bundles on $V$.
A vector bundle $E$ on $V$ is said to be semistable with respect to
$H_1, \ldots, H_{d-1}$ if, for any non-zero subsheaves $G$ of $E$,
\[
\frac{(c_1(G) \cdot c_1(H_1) \cdots c_1(H_{d-1}))}{\rank G} \leq
\frac{(c_1(E) \cdot c_1(H_1) \cdots c_1(H_{d-1}))}{\rank E}.
\]

\bigskip
Let $f : X \to Y$ be a surjective and
projective morphism of quasi-projective varieties
over $k$ with $\dim f = d \geq 1$. 
Let $H_1, \ldots, H_{d-1}$ be line bundles on $X$ and
$E$ a vector bundle on $X$ of rank $r$.
Then, $\left( (2r c_2(E) - (r-1)c_1(E)^2) \cdot c_1(H_1) \cdots c_1(H_{d-1}) \right) \cap [X]$
is a cycle of dimension $\dim Y - 1$ on $X$. Thus, 
$f_*\left( \left( (2r c_2(E) - (r-1)c_1(E)^2) \cdot 
c_1(H_1) \cdots c_1(H_{d-1}) \right) \cap [X]\right)$,
denoted by $\dis_{X/Y}(E; H_1, \ldots, H_{d-1})$, is a divisor on $Y$.
Then, we have the following theorem.

\begin{Theorem}[$\ch(k) = 0$]
\label{thm:nef:psudo:dis:higher}
We assume that $Y$ is smooth over $k$ and 
$H_1, \ldots, H_{d-1}$ are ample.
If $y$ is a point of $Y$,
$f$ is smooth over $y$, and
$E_{\bar{y}}$ is semistable with respect to
$(H_1)_{\bar{y}}, \ldots, (H_{d-1})_{\bar{y}}$
on each connected component of the geometric fiber $X_{\bar{y}}$ over $y$, 
then the discriminant divisor
$\dis_{X/Y}(E; H_1, \ldots, H_{d-1})$ is weakly positive at $y$.
\end{Theorem}

\Proof
We prove this theorem by induction on $d$.
If $d = 1$, then our assertion is nothing more than Corollary~\ref{cor:nef:psudo:dis}.
So we assume $d \geq 2$. 
We choose a sufficiently large integer $n$ so that
$H_{d-1}^{\otimes n}$ is very ample, i.e., there is an embedding
$\iota : X \hookrightarrow \PP^N$ with $H_{d-1}^{\otimes n} \simeq \iota^*(\OO_{\PP^N}(1))$.
By Bertini's theorem, we can find 
a general member $\Gamma \in |\OO_{\PP^N}(1)|$ such that
$X \cap \Gamma$ is integral and $f^{-1}(y) \cap \Gamma$ is smooth.
We set $Z = X \cap \Gamma$ and $g = \rest{f}{Z} :  Z \to Y$.
Since $n$ is sufficiently large,
$g^{-1}(y) \in \left|\rest{H_{d-1}^{\otimes n}}{f^{-1}(y)}\right|$ and
$g^{-1}(y)$ is smooth, by virtue of \cite[Theorem 3.1]{Mo2},
$\rest{E}{Z_{\bar{y}}}$ is semistable with respect to 
$\rest{H_{1}}{Z_{\bar{y}}}, \ldots, \rest{H_{d-2}}{Z_{\bar{y}}}$
on each connected component of $Z_{\bar{y}}$.
Therefore, by hypothesis of induction,
$\dis_{Z/Y}(\rest{E}{Z}; \rest{H_1}{Z}, \ldots, \rest{H_{d-2}}{Z})$ 
is weakly positive at $y$.
On the other hand, we have
\begin{align*}
\dis_{Z/Y}(\rest{E}{Z}; \rest{H_1}{Z}, \ldots, \rest{H_{d-2}}{Z}) & = 
\dis_{X/Y}(E; H_1, \ldots, H_{d-2}, H_{d-1}^{\otimes n}) \\
& = n \dis_{X/Y}(E; H_1, \ldots, H_{d-1}).
\end{align*}
Hence, $\dis_{X/Y}(E; H_1, \ldots, H_{d-1})$ is weakly positive at $y$.
\QED

\section{Relative Bogomolov's inequality in positive characteristic}
\label{sec:relative:bogomolov:inequality:positive:characteristic}
In this section, we will consider a similar result of
Corollary~\ref{cor:nef:psudo:dis} in positive characteristic.
The crucial point of the proof of Corollary~\ref{cor:nef:psudo:dis} is
the semistability of tensor products of semistable vector bundles,
which was studied by a lot of authors \cite{Gi1}, \cite{Ha1}, \cite{Mi}, \cite{Ma1} and etc.
This however does not hold in positive characteristic \cite{Gi0}, so that we will introduce
the strong semistability of vector bundles.

\medskip
Let $C$ be a smooth projective curve over $k$.
For a vector bundle $F$ on $C$, we set
$\mu(F) = \deg(F)/ \rank (F)$, which is
called the {\em slope of $F$}.
A vector bundle $E$ on $C$ is said to be {\em semistable} (resp. {\em stable}) if,
for any proper subbundles $F$ of $E$, $\mu(F) \leq \mu(E)$
(resp. $\mu(F) < \mu(E)$).
Moreover, $E$ is said to be {\em strongly semistable} if,
for any finite morphisms $\pi : C' \to C$
of smooth projective curves over $k$, $\pi^*(E)$ is semistable.
Then, we have the following 
elementary properties of semistable or strongly semistable vector
bundles.

\begin{Proposition}[$\ch(k) \geq 0$]
\label{prop:elem:prop:semistable}
Let $E$ be a vector bundle of rank $r$ on $C$.
\begin{enumerate}
\renewcommand{\labelenumi}{(\arabic{enumi})}
\item
\label{enum:prop:elem:prop:semistable:separable}
Let $\pi : C' \to C$ be a finite separable morphism 
of smooth projective curves over $k$.
If $E$ is semistable, then so is $\pi^*(E)$.

\item
\label{enum:prop:elem:prop:semistable:char:zero}
Under the assumption of $\ch(k) = 0$, 
$E$ is semistable if and only if $E$ is strongly semistable.

\item
\label{enum:prop:elem:prop:semistable:nef}
Let $f : P = \Proj\left( \bigoplus_{m=0}^{\infty} \Sym^m(E) \right) \to C$ 
be the projective bundle of $E$ and
$\OO_P(1)$ the tautological line bundle on $P$. Then,
$E$ is strongly semistable if and only if $\omega_{P/C}^{\otimes -1} = \OO_P(r) \otimes f^*(\det E)^{\otimes -1}$ is
numerically effective.
\end{enumerate}
\end{Proposition}

\Proof
(1) is nothing more than \cite[Lemma~1.1]{Gi1} and
(2) is a consequence of (1).

\medskip
(3) First we assume that $E$ is strongly semistable.
Let $Z$ be any irreducible curves on $P$. If $Z$ is contained in a fiber,
then obviously $(\omega_{P/C}^{\otimes -1} \cdot Z) > 0$. So we may assume that
$Z$ is not contained in any fibers.
Let $C'$ be the normalization of $Z$ and $\pi : C' \to Z \to C$ the induced morphism.
Let $E' = \pi^*(E)$, $f' : P' = \Proj\left( \bigoplus_{m=0}^{\infty} \Sym^m(E') \right) \to C'$ 
the projective bundle of $E'$,
and $\OO_{P'}(1)$ the tautological line bundle on $P'$. Then we have the following commutative
diagram.
\[
\begin{CD}
P @<{\pi'}<< P' \\
@V{f}VV @VV{f'}V \\
C @<<{\pi}< C'
\end{CD}
\]
By our construction, there is a section $Z'$ of $f'$ such that $\pi'(Z') = Z$.
We set $Q' = \rest{\OO_{P'}(1)}{Z'}$. Then, there is a surjective homomorphism
$E' \to Q'$. Since $E'$ is semistable, we have
$\mu(E') \leq \deg(Q')$, which means that
$(\omega_{P'/C'}^{\otimes -1} \cdot Z') \geq 0$. Here, $\omega_{P'/C'} = {\pi'}^*(\omega_{P/C})$.
Thus, we get $(\omega_{P/C}^{\otimes -1} \cdot Z) \geq 0$.

Conversely, we assume that $\omega_{P/C}^{\otimes -1}$ is numerically effective on $P$.
Let $\pi : C' \to C$ be a finite morphism of smooth projective curves over $k$.
We set $f' : P' \to C'$ and $\pi' : P' \to P$ as before.
Then, $\omega_{P'/C'}^{\otimes -1} = {\pi'}^*(\omega_{P/C}^{\otimes -1})$ is
numerically effective on $P'$.
Let $Q$ be a quotient vector bundle of $E' = \pi^*(E)$ with $s = \rank Q$.
The projective bundle $\Proj\left( \bigoplus_{m=0}^{\infty} \Sym^m(Q) \right) \to C'$
gives a subvariety $V'$ of $P'$
with $\deg(Q) = (\OO_{P'}(1)^s \cdot V')$ and $(\OO_{P'}(1)^{s-1} \cdot F' \cdot V') = 1$,
where $F'$ is a fiber of $f'$.
Since $\omega_{P'/C'}^{\otimes -1}$ is numerically effective,
\[
0 \leq \left( (\OO_{P'}(r) \otimes {f'}^*(\det E')^{-1})^s \cdot V' \right)
= r^{s-1}(r \deg(Q) - s \deg(E')).
\]
Thus, $\mu(E') \leq \mu(Q)$.
\QED

First, let us consider symmetric products of strongly semistable vector bundles.

\begin{Theorem}[$\ch(k) \geq 0$]
\label{thm:semistable:sym}
If $E$ is a strongly semistable vector bundle on $C$,
then so is $\Sym^n(E)$ for all $n \geq 0$.
\end{Theorem}

\Proof
Taking a finite covering of $C$, we may assume that
$\deg(E)$ is divisible by $\rank E$.
Let $\theta$ be a line bundle on $C$ with $\deg(\theta) = 1$.
If we set $E_0 = E \otimes \theta^{\otimes -\frac{\deg(E)}{\rank E}}$,
then $\deg(E_0) = 0$ and
$\Sym^n(E_0) = \Sym^n(E) \otimes \theta^{\otimes -\frac{n\deg(E)}{\rank E}}$.
Thus, to prove our theorem, we may assume $\deg(E) = 0$.

We assume that $\Sym^n(E)$ is not strongly semistable for some $n \geq 2$.
By replacing $C$ by a finite covering of $C$, we may assume that
$\Sym^n(E)$ is not semistable.
Let $f : P = \Proj\left( \bigoplus_{n=0}^{\infty} \Sym^n(E) \right) \to C$
be a projective bundle of $E$ and
$\OO_{P}(1)$ the tautological line bundle on $P$.
Let $F$ be the maximal destabilizing sheaf of $\Sym^n(E)$.
In particular, $F$ is semistable and $\mu(F) > 0$. We consider
a composition of homomorphisms
\[
  \alpha : f^*(F) \to f^*(\Sym^n(E)) \to \OO_{P}(n).
\]
Since $f_*(\alpha)$ induces the inclusion $F \to \Sym^n(E)$,
$\alpha$ is a non-trivial homomorphism.

Fix an ample line bundle $A$ on $C$.
Let $l$ be a positive integer with 
$l\mu(F) > n(r-1)\deg(A)$ and $(l, p) = 1$, where $p = \ch(k)$.
Here we claim that $\OO_P(l) \otimes f^*(A)$ is ample.
Let $V$ be an $s$-dimensional subvariety of $P$.
By virtue of Nakai's criterion,
it is sufficient to show $(c_1(\OO_P(l) \otimes f^*(A))^s \cdot V) > 0$.
If $V$ is contained in a fiber, our assertion is trivial. So we may assume that
$V$ is not contained in any fibers.
Then,
\[
(c_1(\OO_P(l) \otimes f^*(A))^s \cdot V) = 
l^s(c_1(\OO_P(1))^s \cdot V) + s l^{s-1} (c_1(\OO_P(1))^{s-1} \cdot c_1(f^*(A)) \cdot V).
\]
Since
$\OO_P(1)$ is numerically effective on $P$ by 
(\ref{enum:prop:elem:prop:semistable:nef}) of
Proposition~\ref{prop:elem:prop:semistable},
$(c_1(\OO_P(1))^s \cdot V) \geq 0$.
Moreover, if $x$ is a general point of $C$,
\[
(c_1(\OO_P(1))^{s-1} \cdot c_1(f^*(A)) \cdot V) =
\deg(A)\deg(\rest{V}{f^{-1}(x)}) > 0.
\]
Therefore, we get our claim.

Thus, there is a positive integer $m$ such that
$(\OO_P(l) \otimes f^*(A))^{\otimes m}$ is very ample and $(m, p) = 1$.
Take general elements $D_1 , \ldots , D_{r-1}$ of
$\left| (\OO_P(l) \otimes f^*(A))^{\otimes m} \right|$
such that $\Gamma = D_1 \cap \cdots \cap D_{r-1}$ is a non-singular curve and
$\rest{f^*(F)}{\Gamma} \to \rest{\OO_{P}(n)}{\Gamma}$ 
is generically surjective.
If $x$ is a general point of $C$,
\iftwelvept
\[
\deg(\Gamma \to C) = (D_1 \cdots D_{r-1} \cdot f^{-1}(x))
= m^{r-1} (c_1(\OO_P(l) \otimes f^*(A))^{r-1} \cdot f^{-1}(x)) =
(ml)^{r-1}.
\]
\else
\begin{align*}
\deg(\Gamma \to C) & = (D_1 \cdots D_{r-1} \cdot f^{-1}(x)) \\
& = m^{r-1} (c_1(\OO_P(l) \otimes f^*(A))^{r-1} \cdot f^{-1}(x)) =
(ml)^{r-1}.
\end{align*}
\fi
Thus, $k(\Gamma)$ is separable over $k(C)$ because $(p, (ml)^{r-1}) = 1$.
Hence, $\rest{f^*(F)}{\Gamma}$ is semistable by
(\ref{enum:prop:elem:prop:semistable:separable}) of
Proposition~\ref{prop:elem:prop:semistable}.
Therefore,
\[
  \frac{(c_1(f^*(F)) \cdot \Gamma)}{\rank F} \leq 
  (c_{1}(\OO_{P}(n)) \cdot \Gamma), 
\]
which implies
\[
  \frac{(c_1(f^*(F)) \cdot c_1(\OO_P(l) \otimes f^*(A))^{r-1})}{\rank F} \leq 
  (c_{1}(\OO_{P}(n)) \cdot c_1(\OO_P(l) \otimes f^*(A))^{r-1}).
\]
This gives rise to
\[
l^{r-1}\mu(F) \leq n(r-1)l^{r-2}\deg(A),
\]
which contradicts to the choice of $l$ with
$l\mu(F) > n(r-1)\deg(A)$.
\QED

As a corollary of Theorem~\ref{thm:semistable:sym}, 
we have the following.

\begin{Corollary}[$\ch(k) \geq 0$]
\label{cor:semistable:tensor}
If $E$ and $F$ are strongly semistable vector bundles on $C$,
then so is $E \otimes F$.
\end{Corollary}

\Proof
Considering a finite covering of $C$ and tensoring line bundles, 
we may assume that
$\deg(E) = \deg(F) = 0$ as in the same way of
the beginning part of the proof of
Theorem~\ref{thm:semistable:sym}.
Then, $E \oplus F$ is strongly semistable.
Thus, by Theorem~\ref{thm:semistable:sym},
$\Sym^2(E \oplus F)$ is strongly semistable.
Here,
\[
\Sym^2(E \oplus F) = (E \otimes F) \oplus \Sym^2(E) \oplus
\Sym^2(F).
\]
Therefore, we can see that $E \otimes F$ is strongly semistable.
\QED

Thus, in the same way as the proof of Corollary~\ref{cor:nef:psudo:dis},
we have the following.

\begin{Corollary}[$\ch(k) \geq 0$]
\label{cor:nef:psudo:dis:in:p}
Let $X$ be a quasi-projective variety over $k$, $Y$ a smooth quasi-projective variety over $k$,
and $f : X \to Y$ a surjective and projective morphism over $k$ with $\dim f = 1$.
Let $E$ be a locally free sheaf on $X$ and
$y$ a point of $Y$.
If $f$ is flat over $y$, 
the geometric fiber $X_{\bar{y}}$ over $y$ is reduced and Gorenstein,
and $E$ is strongly semistable on each
connected component of the normalization of
$X_{\bar{y}}$, 
then $\dis_{X/Y}(E)$ is weakly positive at $y$.
\end{Corollary}

\renewcommand{\thesection}{Appendix \Alph{section}}
\renewcommand{\theTheorem}{\Alph{section}.\arabic{Theorem}}
\renewcommand{\theClaim}{\Alph{section}.\arabic{Theorem}.\arabic{Claim}}
\renewcommand{\theequation}{\Alph{section}.\arabic{Theorem}.\arabic{Claim}}
\setcounter{section}{0}

\section{A certain fibration of hyperelliptic curves}
\label{sec:const:fib:hyperelliptic}

In this section, we would like to construct
a certain fibration of hyperelliptic curves, which is
needed in \S\ref{sec:cone:positive:divisor:moduli:spacc:stable:curve}.
Throughout this section, we assume that $\ch(k) = 0$.

Let us begin with the following lemma.

\begin{Lemma}
\label{lem:conic:fibration}
For non-negative integers $a_1$ and $a_2$,
there are a morphism $f_1 : X_1 \to Y_1$ of smooth projective varieties over $k$,
an effective divisor $D_1$ on $X_1$, a line bundle $L_1$ on $X_1$, and
a line bundle $M_1$ on $Y_1$ with the following properties.
\begin{enumerate}
\renewcommand{\labelenumi}{(\arabic{enumi})}
\item
$\dim X_1 = 2$ and $\dim Y_1 = 1$.

\item
Let $\Sigma_1$ be the set of all critical values of $f_1$, i.e.,
$P \in \Sigma_1$ if and only if $f_1^{-1}(P)$ is a singular variety.
Then, for any $P \in Y_1 \setminus \Sigma_1$, $f_1^{-1}(P)$ is a smooth
rational curve.

\item
$\Sigma_1 \not= \emptyset$, and for any $P \in \Sigma_1$, $f_1^{-1}(P)$ is a reduced curve
consisting of two smooth rational curves $E_P^{1}$ and $E_P^{2}$ joined
at one point transversally.

\item
$D_1$ is smooth over $k$ and
$\rest{f_1}{D_1} : D_1 \to Y_1$ is \'{e}tale.

\item
$(E_P^{1} \cdot D_1) = a_1 + 1$ and $(E_P^{2} \cdot D_1) = a_2 + 1$ for any $P \in \Sigma_1$.
Moreover, $D_1$ does not pass through $E_P^1 \cap E_P^2$.

\item
There is a map $\sigma : \Sigma_1 \to \{ 1, 2 \}$ such that
\[
D_1 \in \left|
L_1^{\otimes a_1 + a_2 + 2} \otimes f_1^*(M_1) \otimes 
\OO_{X_1}\left(-\sum_{P \in \Sigma_1} (a_{\sigma(P)} + 1) E_P^{\sigma(P)}\right) \right|.
\]

\item
$\deg(M_1)$ is divisible by $(a_1 + 1)(a_2 + 1)$.
\end{enumerate}
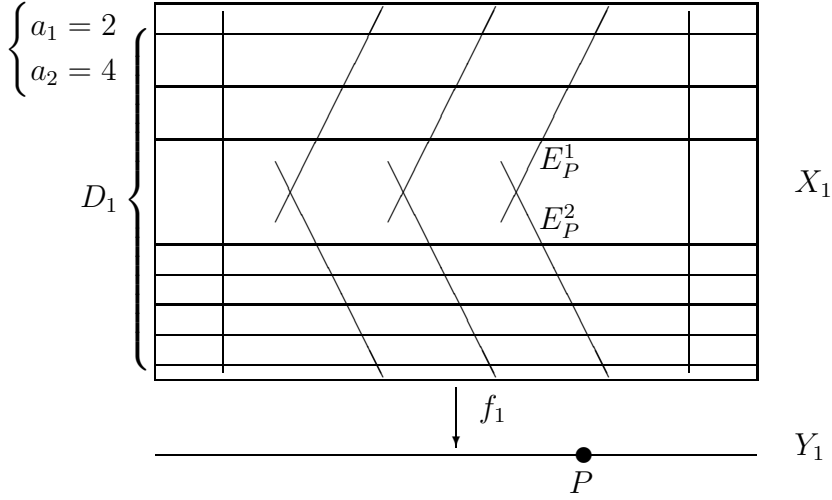
\ifpicture
\bigskip
\begin{center}
\setlength{\unitlength}{1mm}
\begin{picture}(110,70)
\put(10,20){\framebox(80,50){}} \put(95,45){$X_1$}
\put(10,10){\line(1,0){80}} \put(95,10){$Y_1$}
\put(50,19){\vector(0,-1){8}} \put(53,15){$f_1$}
\put(19,21){\line(0,1){48}}
\put(81,21){\line(0,1){48}}
\put(41,41){\line(1,2){14.3}}
\put(41,49){\line(1,-2){14.3}}
\put(26,41){\line(1,2){14.3}} 
\put(26,49){\line(1,-2){14.3}} 
\put(56,41){\line(1,2){14.3}} \put(61,48){$E_P^1$}
\put(56,49){\line(1,-2){14.3}} \put(61,40){$E_P^2$}
\put(10,22){\line(1,0){80}}
\put(10,26){\line(1,0){80}}
\put(10,30){\line(1,0){80}}
\put(10,34){\line(1,0){80}}
\put(10,38){\line(1,0){80}}
\put(10,66){\line(1,0){80}}
\put(10,59){\line(1,0){80}}
\put(10,52){\line(1,0){80}}
\put(0,43){$D_1 \begin{cases} \\ \\ \\ \\ \\ \\ \\ \\ \end{cases}$}
\put(-10,63){$\begin{cases} a_1=2 \\ a_2=4 \end{cases}$}
\put(67,10){\circle*{2}} \put(65,5){$P$}
\end{picture}
\end{center}
\else\fi
\end{Lemma}

\Proof
First of all, let us consider the function $\theta(x)$ defined by
\[
\theta(x) = (a_1 + a_2 + 1) \int_{0}^{x} t^{a_1}(t-1)^{a_2} dt.
\]
Then, $\theta(x)$ is a monic polynomial of degree $a_1 + a_2 + 1$ over $\QQ$.
Moreover, it is easy to see that
\[
\theta'(x) = (a_1 + a_2 + 1) x^{a_1}(x-1)^{a_2},\qquad
\theta(0) = 0\quad\text{and}\quad
\theta(1) = (-1)^{a_2}(a_1+a_2+1)\frac{(a_1)!(a_2)!}{(a_1 + a_2)!}.
\]
Thus, there are distinct non-zero algebraic numbers
$\alpha_1, \ldots, \alpha_{a_2}$ and $\beta_1, \ldots, \beta_{a_1}$ such that
\[
\theta(x) = x^{a_1+1}(x-\alpha_1) \cdots (x-\alpha_{a_2})
\]
and
\[
\theta(x) - \theta(1) = (x-1)^{a_2+1}(x-1-\beta_1) \cdots (x-1-\beta_{a_1}).
\]
Here we set
\[
F(X,Y) = Y^{a_1 + a_2 + 1}\theta(X/Y) = X^{a_1 + 1}(X - \alpha_1 Y) \cdots (X - \alpha_{a_2} Y)
\]
and
\[
G(X,Y,S,T) = T F(X,Y) - Y^{a_1 + a_2 + 1} S.
\]
Then, $F$ is a homogeneous polynomial of degree $a_1 + a_2 + 1$ over $\QQ$, and
$G$ is a bi-homogeneous polynomial of bi-degree $(a_1+ a_2 + 1, 1)$ in 
$\QQ[X,Y] \otimes_{\QQ} \QQ[S,T]$.
Let $D'$ (resp. $D''$) be the curve on $\PP^1_{(X,Y)} \times \PP^1_{(S,T)}$ 
given by the equation $\{ G = 0 \}$ (resp. $\{ Y = 0 \}$),
where $\PP^1_{(X,Y)} = \Proj(k[X, Y])$ and $\PP^1_{(S,T)} = \Proj(k[S, T])$.
Moreover, we set $D = D' + D''$.
Let $p : \PP^1_{(X,Y)} \times \PP^1_{(S,T)} \to \PP^1_{(X,Y)}$ and
$q : \PP^1_{(X,Y)} \times \PP^1_{(S,T)} \to \PP^1_{(S,T)}$
be the natural projections. 
Then, $D'$ (resp. $D''$) is an element of 
the linear system $\left| p^*(\OO_{\PP^1}(a_1 + a_2 + 1)) \otimes q^*(\OO_{\PP^1}(1)) \right|$
(resp $\left| p^*(\OO_{\PP^1}(1)) \right|$). Thus, 
$D \in \left| p^*(\OO_{\PP^1}(a_1 + a_2 + 2)) \otimes q^*(\OO_{\PP^1}(1)) \right|$,
$(D' \cdot D'') = 1$
and $D' \cap D'' = \{ ((1 : 0), (1 : 0)) \}$.
Here we claim the following.

\begin{Claim}
\label{claim:property:of:D:prime}
\begin{enumerate}
\renewcommand{\labelenumi}{(\alph{enumi})}
\item
$D'$ is a smooth rational curve.

\item
Let $\pi' : D' \to \PP^1_{(S,T)}$ be the morphism induced by
the projection $q : \PP^1_{(X,Y)} \times \PP^1_{(S,T)} \to \PP^1_{(S,T)}$.
If we set 
$Q_1 = ((0 : 1), (0 : 1))$, $Q_2 = ((1 : 1), (\theta(1), 1))$ and
$Q_3 = ((1 : 0), (1 : 0))$, then
the set of ramification points of $\pi'$ is $\{ Q_1, Q_2, Q_3 \}$.
Further, the ramification indexes at $Q_1$, $Q_2$ and $Q_3$ are
$a_1 + 1$, $a_2 + 1$ and $a_1 + a_2 + 1$ respectively.
\end{enumerate}
\end{Claim}

\Proof
(a) Since $F(X, Y)$ has no factor of $Y$, the morphism
$e : \PP^1_{(X,Y)} \to D'$ given by
\[
e(x:y) = \left( (x : y), (F(x,y) : y^{a_1+a_2+1}) \right)
\]
is well defined.
Moreover, if we set $e' = \rest{p}{D'}$, then it is easy to see that
$e \cdot e' = \operatorname{id}_{D'}$ and $e' \cdot e = \operatorname{id}_{\PP^1}$.
Thus, $D'$ is a smooth rational curve.

\medskip
(b) Pick up a point $(\lambda : \mu) \in \PP^1_{(S, T)}$.
Then,
$G_{(\lambda, \mu)} = \mu F(X, Y) - Y^{a_1+a_2+1} \lambda$ is a homogeneous polynomial of
degree $a_1 + a_2 +1$. 

First, we assume that $\mu \not= 0$, hence we may assume that $\mu=1$.
Then, $Y$ is not a factor of $G_{(\lambda, 1)}(X, Y)$, which means that
${\pi'}^{-1}((\lambda : 1))$ sits in the affine open set
$\Spec(k[X/Y,S/T])$. Thus,
\[
{\pi'}^{-1}((\lambda : 1)) = \{ ((\gamma : 1), (\lambda : 1) ) \mid \theta(\gamma) - \lambda = 0 \}.
\]
Hence, in order to get ramification points of $\pi'$, we need to see
multiple roots of $\phi(x) = \theta(x) - \lambda$.
Here we will check that $\phi(x)$ has a multiple root
if and only if $\lambda$ is either $0$ or $\theta(1)$.
Moreover, if $\lambda$ is $0$ (resp. $\theta(1)$), then
$0$ (resp. $1$) is the only multiple root of $\phi(x)$ with multiplicity $a_1 + 1$
(resp. $a_2+1$).

Let $\gamma$ be a multiple root of $\phi(x) = 0$.
Then, $\phi(\gamma) = \phi'(\gamma) = 0$. Here,
\[
\phi'(x) = (a_1+a_2+1)x^{a_1}(x-1)^{a_2}.
\]
Thus, $\gamma$ is either $0$ or $1$.
If $\gamma = 0$, then $\lambda = \theta(0) = 0$.
If $\gamma = 1$, then $\lambda = \theta(1)$. 
In the same way, we can easily check the remaining part of our assertion.

Therefore, we get two ramification points $Q_1$ and $Q_2$ whose ramification indexes
are $a_1+1$ and $a_2 +1$  respectively.

Next, we assume that $\mu = 0$, hence we may assume $\lambda=1$.
Then, $G_{(\lambda, \mu)} =  - Y^{a_1+a_2+1}$.  Thus, $P_3$ is a ramification point 
whose ramification index is $a_1+a_2+1$.
\QED

\medskip
Here we set $P_i = q(Q_i)$ ($i=1,2,3$),
$b_1 = a_1+1$, $b_2 = a_2 + 1$, and $b_3 = a_1+a_2+1$.

\begin{Claim}
\label{claim:cyclic:cover}
There is a cyclic covering $h_1 : Y_1 \to \PP^1_{(S,T)}$
of smooth projective curves such that
$\deg(h_1) = b_1 b_2 b_3$ and that,
for any $i=1,2,3$ and any $P \in h_1^{-1}(P_i)$,
the ramification index of $h_1$ at $P$ is $b_i$.
\end{Claim}

\Proof
Since $b_1b_2 + b_2b_3 + b_3b_1 \leq 3 b_1b_2b_3$,
there is an effective and reduced divisor $d$ on $\PP^1_{(S,T)}$
such that $P_i \not\in \Supp(d)$ for each $i=1,2,3$ and
\[
b_2b_3 P_1 + b_3b_1 P_2+ b_1b_2 P_3 + d \in \left| \OO_{\PP^1}(3b_1b_2b_3) \right|.
\]
Let $w$ be a section of $H^0(\OO_{\PP^1}(3b_1b_2b_3))$ with
$\zero(w) = b_2b_3 P_1+ b_3b_1 P_2 + b_1b_2 P_3+ d$.
Then, $w$ gives rise to the ring structure on
$\bigoplus_{i=0}^{b_1b_2b_3 - 1} \OO_{\PP^1}(-3i)$.
Let $Y_1$ be the normalization of 
\[
\Spec\left( \bigoplus_{i=0}^{b_1b_2b_3 - 1} \OO_{\PP^1}(-3i) \right)
\]
and $h_1: Y_1 \to \PP^1$ the induced morphism.
Then, by our choice of $w$, it is easy to see that
$h_1 : Y_1 \to \PP^1$ satisfies the desired properties.
\QED

Let $p_1 : \PP^1_{(X,Y)} \times Y_1 \to \PP^1_{(X,Y)}$ and
$q_1 : \PP^1_{(X,Y)} \times Y_1 \to Y_1$
be the natural projections, and 
$u_1 = \operatorname{id} \times h_1 : 
\PP^1_{(X,Y)} \times Y_1 \to  \PP^1_{(X,Y)} \times \PP^1_{(S,T)}$. 
Then, we have a commutative diagram:
\[
\begin{CD}
\PP^1_{(X,Y)} \times \PP^1_{(S,T)} @<{u_1}<< \PP^1_{(X,Y)} \times Y_1 \\
@V{q}VV @VV{q_1}V \\
\PP^1_{(S,T)} @<<{h_1}< Y_1
\end{CD}
\]
We set $h_1^{-1}(P_1)$, $h_1^{-1}(P_2)$ and $h_1^{-1}(P_3)$ as follows.
\[
\begin{cases}
h_1^{-1}(P_1) = \{ P_{1,1}, \ldots, P_{1,b_2b_3} \}, \\
h_1^{-1}(P_2) = \{ P_{2,1}, \ldots, P_{1,b_3b_1} \},  \\
h_1^{-1}(P_3) = \{ P_{3,1}, \ldots, P_{3,b_1b_2} \}.
\end{cases}
\]
Then, there is a unique $Q_{i,j}$ on $\PP^1_{(X,Y)} \times Y_1$
with $q_1(Q_{i,j}) = P_{i,j}$ and $u_1(Q_{i,j}) = Q_i$.

\begin{Claim}
\label{claim:sing:of:D}
\begin{enumerate}
\renewcommand{\labelenumi}{(\alph{enumi})}
\item
$u_1^*(D)$ is \'{e}tale over $Y_1$ outside $\{ Q_{i,j} \}_{i,j}$.
In particular, $u_1^*(D)$ is smooth over $k$ outside $\{ Q_{i,j} \}_{i,j}$.

\item
If we set $c_1 = a_1+ 1$, $c_2 = a_2 + 1$ and $c_3 = a_1 + a_2 + 2$,
then $u_1^*(D)$ has an ordinary $c_i$-fold point at $Q_{i,j}$ for every
$i,j$. Moreover, each tangent of $u_1^*(D)$ at $Q_{i,j}$ is different from
the fiber $q_1^{-1}(P_{i,j})$.
\end{enumerate}
\end{Claim}

\Proof
(a) is trivial because $\rest{q}{D} : D \to \PP^1_{(S,T)}$ is \'{e}tale outside
$\{ Q_1, Q_2, Q_3\}$.
Since $u_1^*(D'') = p_1^{-1}((1:0))$, in order to see (b),
it is sufficient to check the following.
$u_1^*(D')$ has an ordinary $b_i$-fold point at $Q_{i,j}$ for every
$i,j$. Moreover, for i=1,2, each tangent of $u_1^*(D')$ at $Q_{i,j}$ is different from
the fiber $q_1^{-1}(P_{i,j})$, and each tangent of $u_1^*(D')$ at $Q_{3,j}$ is different from
the fiber $q_1^{-1}(P_{3,j})$ and $p_1^{-1}((1:0))$.

First we assume $i=1$. Let $z$ be a local parameter of $Y_1$ at $P_{1,j}$,
$x = X/Y$, and $s = S/T$.
Then, $(x,z)$ gives a local parameter of $\PP^1_{(X,Y)} \times Y_1$ at $Q_{1,j}$.
Since $s = v(z) z^{a_1+1}$ for some $v(z)$ with $v(0) \not= 0$,
$u_1^*(D')$ is defined by
\[
x^{a_1+1}(x - \alpha_1) \cdots (x - \alpha_{a_2}) - v(z) z^{a_1+1} = 0
\]
around $Q_{1, j}$. Thus, since $\alpha_1 \cdots \alpha_{a_2} \not= 0$,
$Q_{1,j}$ is an ordinary $(a_1+1)$-fold point and
each tangent is different from $\{ z = 0 \}$.

Next we assume $i=2$. Let $z$ be a local parameter of $Y_1$ at $P_{2,j}$,
$x' = X/Y - 1$, and $s' = S/T - \theta(1)$.
Then, $(x',z)$ gives a local parameter of $\PP^1_{(X,Y)} \times Y_1$ at $Q_{2,j}$.
Since $s' = v(z) z^{a_2+1}$ for some $v(z)$ with $v(0) \not= 0$,
$u_1^*(D')$ is defined by
\[
(x')^{a_2+1}(x' - \beta_1) \cdots (x' - \beta_{a_1}) - v(z) z^{a_2+1} = 0
\]
around $Q_{2, j}$. Thus, we can see our assertion in this case
because $\beta_1 \cdots \beta_{a_1} \not= 0$.

Finally we assume that $i=3$.
Let $z$ be a local parameter of $Y_1$ at $P_{3,j}$,
$y = Y/X$, and $t = T/S$.
Since $t = v(z) z^{a_1+a_2+1}$ for some $v(z)$ with $v(0) \not= 0$,
$u_1^*(D')$ is defined by
\[
v(z) z^{a_1+a_2+1} (1 - \alpha_1 y) \cdots (1 - \alpha_{a_2} y) - y^{a_1+a_2+1} = 0
\]
around $Q_{3, j}$.
Thus, $Q_{3,j}$ is an ordinary $(a_1+a_2+1)$-fold point and
each tangent is different from $\{ z = 0 \}$ and $\{ y = 0 \}$.
\QED

Let $\mu_1 : Z_1 \to \PP^1_{(X,Y)} \times Y_1$ be blowing-ups at all points $Q_{i,j}$, and
$E_{i,j}$ $(-1)$-curve over $Q_{i,j}$.
Let $\overline{D}_1$ be the strict transform of $u_1^*(D)$ by $\mu_1$, and
$g_1 = q_1 \cdot \mu_1$. Then,
by the previous claim, 
$\overline{D}_1$ is \'{e}tale over $Y_1$ and
\[
\overline{D}_1 \in \left| 
\mu_1^*(p_1^*(\OO_{\PP^1}(1)))^{\otimes a_1+a_2+2} \otimes 
g_1^*(h_1^*(\OO_{\PP^1}(1))) \otimes 
\OO_{Z_1}\left(-\sum_{i,j} c_i E_{i,j}\right) \right|.
\]
Let $F_j$ be the strict transform of the fiber $q_1^{-1}(P_{3,j})$.
Note that $F_j \cap \overline{D}_1 = \emptyset$ for all $j$.
Since $F_j$'s are $(-1)$-curve, we can contract them.
Let $\nu_1 : Z_1 \to X_1$ be the contraction of $F_j$'s, and
$f_1 : X_1 \to Y_1$ the induced morphism.
\[
\begin{CD}
\PP^1_{(X,Y)} \times \PP^1_{(S,T)} @<{u_1}<< 
\PP^1_{(X,Y)} \times Y_1 @<{\mu_1}<< Z_1 @>{\nu_1}>> X_1 \\
@V{q}VV @V{q_1}VV @V{g_1}VV @V{f_1}VV \\
\PP^1_{(S,T)} @<<{h_1}< Y_1 @= Y_1 @= Y_1
\end{CD}
\]
Here we set $D_1 = (\nu_1)_*(\overline{D}_1)$,
$L_1 = (\nu_1)_*(\mu_1^*(p_1^*(\OO_{\PP^1}(1))))^{**}$, and
\[
M_1 = h_1^*(\OO_{\PP^1}(1)) \otimes \OO_{Y_1}\left(-\sum_j c_3 P_{3,j}\right) \simeq
\OO_{Y_1}\left(-\sum_j P_{3,j}\right).
\]
Then, since
$(\nu_1)_*(g_1^*(h_1^*(\OO_{\PP^1}(1)))) =
f_1^*(h_1^*(\OO_{\PP^1}(1)))$ and $(\nu_1)_*(E_{3,j}) = f_1^*(P_{3,j})$
for all $j$, we can see
\[
D_1 \in \left| L_1^{\otimes a_1+a_2+2} \otimes 
f_1^*(M_1) \otimes 
\OO_{X_1}\left(-\sum_{\substack{i=1,2, \\ j \geq 1}} b_i \nu_1(E_{i,j})\right) \right|.
\]
Therefore, by our construction of $f_1 : X_1 \to Y_1$, $D_1$, $L_1$ and $M_1$,
it is easy to see all properties (1) --- (6) in Lemma~\ref{lem:conic:fibration}.
\QED

\begin{Proposition}
\label{prop:hyperelliptic:fibration}
Let $g$ and $a$ be integers with $g \geq 1$ and $0 \leq a \leq [g/2]$.
Then, there are a smooth projective surface $X$ over $k$,
a smooth projective curve $C$ over $k$, and a surjective morphism
$f : X \to Y$ over $k$ with the following properties.
\begin{enumerate}
\renewcommand{\labelenumi}{(\arabic{enumi})}
\item
The generic fiber of $f$ is a smooth hyperelliptic curve of genus $g$.

\item
$f$ is not smooth and every fiber is reduced.
Moreover, every singular fiber of $f$ is
a nodal curve consisting of a smooth curve of genus $a$ and
a smooth curve of genus $g-a$ joined at one point.
\end{enumerate}
\end{Proposition}

\Proof
Applying Lemma~\ref{lem:conic:fibration} to the case where $a_1= 2a$ and $a_2 = 2g-2a$,
we fix a conic fibration as in Lemma~\ref{lem:conic:fibration}.
Adding one point to $\Sigma_1$, if necessarily, we can take an effective and reduced
divisor $d$ on $Y_1$ such that $\Sigma_1 \subseteq \Supp(d)$ and $\deg(d)$ is even.
Thus, there is a line bundle $\vartheta$ on $Y_1$ with $\OO_{Y_1}(d) \simeq \vartheta^{\otimes 2}$,
which produces a double covering $h_2 : Y \to Y_1$ of smooth projective curves such that
$h_2$ is ramified over $\Sigma_1$.
Let $\mu_2 : X_2 \to X_1 \times_{Y_1} Y$ be the minimal resolution of singularities of 
$X_1 \times_{Y_1} Y$.
We set the induced morphisms as follows.
\[
\begin{CD}
X_1 @<{u_2}<< X_2 \\
@V{f_1}VV @VV{f_2}V \\
Y_1 @<{h_2}<< Y
\end{CD}
\]
Let $\Sigma_2$ be the set of all critical values of $f_2$. 
Here, for all $Q \in \Sigma_2$,
$f_2^{-1}(Q)$ is reduced, and there is the irreducible decomposition
$f_2^{-1}(Q) = \overline{E}_Q^1 + \overline{E}_Q^2 + B_Q$ such that
$u_2(\overline{E}_Q^i) = E_{h_2(Q)}^i$ for $i=1, 2$ and
$B_Q$ is a $(-2)$-curves.
We set $D_2 = u_2^*(D_1)$ and $B = \sum_{Q \in \Sigma_2} B_Q$. 
Then, $D_2$ is \'{e}tale over $Y$ and
$D_2+ B$ is smooth over $k$ because $D_2 \cap B = \emptyset$.
Moreover, 
\[
D_2 \in \left|
u_2^*(L_1)^{\otimes 2g+2} \otimes f_2^*(h_2^*(M_1)) \otimes 
\OO_{X_2}\left(-u_2^*\left(\sum_{P \in \Sigma_1} (a_{\sigma(P)} + 1) E_P^{\sigma(P)}
\right) \right) \right|.
\]
Let $\sigma_2 : \Sigma_2 \to \{ 1, 2 \}$ be the map given by
$\sigma_2 = \sigma \cdot h_2$.
Then
\[
u_2^*\left(\sum_{P \in \Sigma_1} (a_{\sigma(P)} + 1) E_P^{\sigma(P)} \right) =
\sum_{Q \in \Sigma_2} (a_{\sigma_2(Q)} + 1) (2 \overline{E}_Q^{\sigma_2(Q)} +  B_Q).
\]
Therefore,
\[
D_2 + B \in \left|
u_2^*(L_1)^{\otimes 2g+2} \otimes f_2^*(h_2^*(M_1)) \otimes 
\OO_{X_2}\left(- \sum_{Q \in \Sigma_2} ( 2(a_{\sigma_2(Q)} + 1) \overline{E}_Q^{\sigma_2(Q)} +
a_{\sigma_2(Q)} B_Q ) \right) \right|.
\]
Here, since $\deg(h_2^*(M_1)) = 2 \deg(M_1)$, $h_2^*(M_1)$ is divisible by $2$ in
$\Pic(Y)$. Further, $a_i$ is even for each $i=1, 2$.
Thus, 
\[
u_2^*(L_1)^{\otimes 2g+2} \otimes f_2^*(h_2^*(M_1)) \otimes 
\OO_{X_2}\left(-  \sum_{Q \in \Sigma_2} ( 2(a_{\sigma_2(Q)} + 1) \overline{E}_Q^{\sigma_2(Q)} +
a_{\sigma_2(Q)} B_Q ) \right)
\]
is divisible by $2$ in $\Pic(X_2)$, i.e., there is a line bundle $H$ on $X_2$ with
\[
H^{\otimes 2} \simeq
u_2^*(L_1)^{\otimes 2g+2} \otimes f_2^*(h_2^*(M_1)) \otimes 
\OO_{X_2}\left(-  \sum_{Q \in \Sigma_2} ( 2(a_{\sigma_2(Q)}+1) \overline{E}_Q^{\sigma_2(Q)} +
a_{\sigma_2(Q)} B_Q ) \right).
\]
Hence, we can construct a double covering $\mu_3 : X_3 \to X_2$ of
smooth projective surfaces such that
$\mu_3$ is ramified over $D_2 + B$.
Let $f_3 : X_3 \to Y$ be the induced morphism.
Then, there is the irreducible decomposition
\[
f^{-1}(Q) = \overline{C}^1_Q + \overline{C}^2_Q + 2 \overline{B}_Q
\]
as cycles such that $\mu_3(\overline{C}^i_Q) = \overline{E}^i_Q$ ($i=1,2$) and
$\mu_3(\overline{B}_Q) = B_Q$.
Here it is easy to check that $\overline{B}_Q$ is
a $(-1)$-curve. Thus, we have the contraction $\nu_3 : X_3 \to X$ of $\overline{B}_Q$'s, and
the induced morphism $f : X \to Y$.
\[
\begin{CD}
X_1 @<{u_2}<< X_2@<{\mu_3}<< X_3 @>{\nu_3}>> X \\
@V{f_1}VV @V{f_2}VV @V{f_3}VV @V{f}VV \\
Y_1 @<{h_2}<< Y @= Y @= Y
\end{CD}
\]
We denote $\nu_3(\overline{C}^i_Q)$ by $C_Q^i$. Then,
$C^1_Q$ (resp. $C^2_Q$) is a smooth projective curve of genus $a$ (resp. $g-a$),
$(C_Q^1 \cdot C_Q^2) = 1$, and $f^{-1}(Q) = C^1_Q + C^2_Q$.
Thus, $f : X \to Y$ is our desired fibration.
\QED

In the same way, we can also show the following proposition.

\begin{Proposition}
\label{prop:hyperelliptic:fibration:2}
Let $g$ and $a$ be integers with $g \geq 1$ and $0 \leq a \leq [(g-1)/2]$.
Then, there are a smooth projective surface $X$ over $k$,
a smooth projective curve $C$ over $k$, and a surjective morphism
$f : X \to Y$ over $k$ with the following properties.
\begin{enumerate}
\renewcommand{\labelenumi}{(\arabic{enumi})}
\item
The generic fiber of $f$ is a smooth hyperelliptic curve of genus $g$.

\item
$f$ is not smooth and every fiber is reduced.
Moreover, every singular fiber of $f$ is
a nodal curve consisting of a smooth curve of genus $a$ and
a smooth curve of genus $g-a-1$ joined at two points.
\end{enumerate}
\end{Proposition}

\Proof
Applying Lemma~\ref{lem:conic:fibration} to the case where $a_1 = 2a+1$ and $a_2 = 2g-2a - 1$,
we fix a conic fibration as in Lemma~\ref{lem:conic:fibration}.
In this case, $\deg(M_1)$ is even. Thus,
\[
L_1^{\otimes 2g+2} \otimes f_1^*(M_1) \otimes 
\OO_{X_1}\left(-\sum_{P \in \Sigma_1} (a_{\sigma(P)} + 1) E_P^{\sigma(P)}\right) 
\]
is divisible by $2$ in $\Pic(X_1)$.
Therefore, there is a double covering $\mu : X \to X_1$ of
smooth projective surfaces such that
$\mu_3$ is ramified over $D_1$.
Then, the induced morphism $f : X \to Y_1$ is a desired fibration.
\QED

\end{document}